\documentclass[12pt]{article}
\title{Spectrum of the Laplacian on Quaternionic K\"ahler Manifolds}

\author{Shengli Kong, Peter Li\thanks{Research partially supported by NSF grant DMS-0503735}  and Detang Zhou\thanks{Research partially supported by CAPES and CNPq of Brazil.}}

\usepackage{amsmath,amsthm}
\usepackage{graphicx}
\usepackage{amsmath}
\usepackage{amsfonts}
\usepackage{ae} % or {zefonts}
\usepackage[T1]{fontenc}
\usepackage[ansinew]{inputenc}
\usepackage{amssymb}

\makeatletter
 \@addtoreset{equation}{section}
\makeatother
\newtheorem{thm}{Theorem}[section]
\newtheorem{lem}{Lemma}[section]
\newtheorem{cor}{Corollary}[section]
\newtheorem{rema}{Remark}[section]

\newtheorem{defn}{Definition}[section]

\newcommand{\End}{\textrm{End}}
\newcommand{\Ric}{\textrm{\rm Ric}}
\newcommand{\K}{\mathcal{K}}
\newcommand{\la}{\langle}
\newcommand{\ra}{\rangle}

\newcommand\p{\partial}
\newcommand\n{\nabla}

\renewcommand{\proof}[1]{\noindent {\bf Proof#1:\  }}
\newcommand{\QED}{\hfill$\Box$\medskip}
\newcommand{\kk}{4(n+2)}
\begin{document}
\date{}
\maketitle

 \begin{abstract}Let $M^{4n}$ be a complete quaternionic K\"ahler manifold
 with scalar curvature bounded below by $-16n(n+2)$. We get a sharp
 estimate for the first eigenvalue $\lambda_1(M)$ of the Laplacian
 which is $\lambda_1(M)\le (2n+1)^2$. If the equality holds, then either
 $M$ has only one end, or $M$ is diffeomorphic to $\mathbb{R}\times N$ with N given
  by a compact manifold.
  Moreover, if $M$ is of bounded curvature,  $M$ is covered by the quaterionic hyperbolic space
  $\mathbb{QH}^n$ and $N$ is a compact quotient of the generalized Heisenberg group.  When  $\lambda_1(M)\ge
 \frac{8(n+2)}3$,  we also prove that $M$ must have only one end with infinite volume.

 \end{abstract}

\baselineskip=15pt

\section* {0 Introduction}\label{intro}
Let $M^n$ be a complete $n$-dimensional Riemannian manifold  whose
Ricci curvature bounded below by $-(n-1)$. It is well  known from
Cheng \cite{Ch} that the first eigenvalue $\lambda_1(M)$ satisfies
\begin{equation*}
    \lambda_1(M)\le \frac{(n-1)^2}4.
\end{equation*}

In \cite{LW3}, Li and Wang proved an analogous theorem for complete K\"ahler manifolds.  They showed that if $M^{2n}$ is a complete
K\"ahler  manifold of complex dimension
$n$ with holomorphic bisectional curvature ${\rm BK}_M$ bounded below by $-1$, then  the first eigenvalue $\lambda_1(M)$
satisfies
\begin{equation*}
    \lambda_1(M)\le n^2.
\end{equation*}
Here ${\rm BK}_M \ge -1$ means that
\[
R_{i\bar i j \bar j} \ge -(1 +\delta_{ij})
\]
for any unitary frame $e_1, \dots, e_n.$

In this paper, we prove the corresponding Laplacian comparison
theorem for a  quaterionic K\"ahler manifold $M^{4n}$. As an
application we get the sharp estimate $\lambda_1(M) $ for a complete
quaterionic K\"ahler manifold $M^{4n}$ with scalar curvature bounded
below by $-16n(n+2)$ as
\begin{equation*}
    \lambda_1(M)\le (2n+1)^2.
\end{equation*}

It is an interesting question to ask what one can say about those
manifolds when the above inequalities are realized as equalities. In
works of Li and Wang \cite{LW1} and \cite{LW2}, the authors obtained
the following theorems.  The first was a generalization of the
theory of Witten-Yau \cite{WY}, Cai-Galloway \cite{CG}, and Wang
\cite{W} for conformally compact manifolds.  The second was to
answer the aforementioned question.

\begin{thm}Let $M^n$ be a complete Riemannian manifold of dimension
$n\ge 3$ with Ricci curvature bounded below by $-(n-1)$. If
$\lambda_1(M)\ge n-2$, then either \newline $(1)$ $M$ has only one
infinite volume end; or \newline $(2)$ $M=\mathbb{R}\times N$ with
warped product metric of the form
\[
ds_M^2=dt^2+\cosh^2t \,ds_N^2,
\]
where $N$ is an $(n-1)$-dimensional compact manifold of Ricci curvature bounded
below by $\lambda_1(M)$.
\end{thm}

\begin{thm}Let $M^n$ be a complete Riemannian manifold of dimension
$n\ge 2$ with Ricci curvature bounded below by $-(n-1)$. If
$\lambda_1(M)\ge \frac{(n-1)^2}4$, then either \newline $(1)$ $M$
has no finite volume end; or \newline $(2)$ $M=\mathbb{R}\times N$
with warped product metric of the form
\[
ds_M^2=dt^2+e^{2t}\, ds_N^2,\]
where $N$ is an $(n-1)$-dimensional compact manifold of nonnegative Ricci
curvature.
\end{thm}

In \cite{LW3} and \cite{LW5}, Li and Wang also consider the K\"ahler
case.  They proved the following theorems.

\begin{thm}Let $M^n$ be a complete K\"ahler manifold of complex dimension
$n\ge 1$ with Ricci curvature bounded below by
\[
\Ric_M \ge -2(n+1).
\]
If
$\lambda_1(M)> \frac{n+1}2$, then  $M$ must have only one infinite volume end.
\end{thm}

\begin{thm}Let $M^n$ be a complete K\"ahler manifold of complex dimension
$n\ge 2$ with holomorphic bisectional curvature bounded by
\[
{\rm BK}_M \ge -1.
\]
If $\lambda_1(M)\ge n^2$, then either \newline $(1)$ $M$ has only
one end; or \newline $(2)$ $M=\mathbb{R}\times N$ with
 $N$ being a compact manifold.  Moreover the metric on $M$ is of the form
 \[
 ds^2_M = dt^2 + e^{4t}\,\omega_2^2 + e^{2t}\,\sum_{i=3}^{2n} \omega_{i}^2,
 \]
 where $\{\omega_2, \omega_3, \dots, \omega_{2n}\}$ are orthonormal coframe of $N$ with $J dt = \omega_2.$  If $M$ has bounded curvature, then we further conclude that $M$ is covered by $\mathbb{CH}^n$ and $N$ is a compact quotient of the Heisenberg group.
\end{thm}

In \cite{LW5}, the authors pointed out that the assumption on the lower bound of $\lambda_1(M)$ in  Theorem 0.3 is sharp, since one can construct  $M$ of the form $M = \Sigma \times N$ satisfying
\begin{eqnarray} \label{eqn0.1}
\Ric_M \ge -2(n+1)
\end{eqnarray}
and
\begin{eqnarray} \label{eqn0.2}
\lambda_1(M) = \frac{n+1}2
\end{eqnarray}
with $N$ being a compact K\"ahler manifold and $\Sigma$ being a
complete surface with at least two infinite volume ends.  However,
it is still an open question to characterized all those complete
K\"ahler manifolds satisfying conditions $(\ref{eqn0.1})$ and
$(\ref{eqn0.2})$.

In sections 4 and 5, we will prove the following quaternionic
K\"ahler versions of the above theorems.
\begin{thm}\label{thmnonp}Let $(M^{4n}, g)$  be a complete quaternionic K\"ahler
manifold with scalar curvature satisfying
\[ {\rm S}_M\ge -16n(n+2).
\]
If $\lambda_1(M)\ge\frac{8(n+2)}{3}$, then $M$ must have only one
infinite volume end.
\end{thm}

\begin{thm}\label{thmend}Let $(M^{4n}, g)$  be a complete quaternionic K\"ahler
manifold with  scalar curvature satisfying
\[ {\rm S}_M\ge -16n(n+2).
\]
If $\lambda_1(M)\ge (2n+1)^2$,  then either \newline $(1)$ $M$ has
only one end, or \newline $(2)$ $M$ is diffeomorphic to
$\mathbb{R}\times N$ where N is a compact manifold.  Moreover, the
metric is given by the form
 \[ ds_M^2 = dt^2 + e^{4t}\,\sum_{p=2}^4 \omega_p^2 + e^{2t}\,\sum_{\alpha=5}^{4n} \omega_{\alpha}^2,
 \]
 where $\{\omega_2, \dots, \omega_{4n}\}$ are orthonormal coframes for $N.$
 If $M$ is of bounded curvature then we further conclude that $M$ is covered by the quaterionic hyperbolic space
  $\mathbb{QH}^n$ and $N$ is a compact quotient of the generalized Heisenberg group.
\end{thm}
\begin{rema}
It is known that a  horosphere in $\mathbb{QH}^n$ is isometric to a
certain generalized Heisenberg group with three-dimensional center
and left-invariant Riemannian metric. Such generalized Heisenberg
groups have compact quotients. For an explicit construction see for
instance Example 2.6 in \cite{G}. We don't have an example to show
that the bounded curvature condition in Theorem \ref{thmend} is
necessary. If such an example exists, its curvature should decay at
exponentially in some directions.
\end{rema}
 Perhaps it is interesting to
restrict our attention to the special case when $M^{4n}=
\mathbb{QH}^n/\Gamma$ is given by the quotient of the quaternionic
hyperbolic space $\mathbb {QH}^n$ with a discrete group of isometies
$\Gamma$.  In particular, it is instructional to compare with
previous results by Corlette \cite{C2} and Corlette-Iozzi \cite{CI}
where  Lie group theoretic approach was used in understanding these
manifolds.  For example, in \cite{CI}, the authors proved a
Patterson-Sullivan type formula for $\lambda_1(M)$ in terms of the
Hausdorff dimension $\delta(\Gamma)$ of the limit set of $\Gamma.$
More specifically, they proved that if $\Gamma$ is geometrically
finite, then for $\delta(\Gamma) \ge 2n+1$ one has
$$
\lambda_1(M) = \delta(\Gamma) ((4n+2) - \delta(\Gamma)).
$$
Hence in this case, the condition in Theorem \ref{thmend} on $\lambda_1(M) = (2n+1)^2$ is equivalent to the condition $\delta(\Gamma) = 2n+1.$

In \cite{C2} (Theorem 4.4), Corlette also pointed out that by a result of Kostant $\lambda_1(M) = 0$ or $\lambda_1(M)\ge 8n$.  On the other hand, it was also shown in \cite{CI} that if $\Gamma$ is geometrically finite and torsion free, then $M = \mathbb {QH}^n/ \Gamma$ must have at most one  end with infinite volume.  These two statements give an interesting comparison to  Theorem \ref{thmnonp} stated above.

We would also like to point out to the interested readers that in \cite{LW4} and \cite{LW5} Li and Wang considered a more general class of manifolds satisfying a weighted Poincar\'e inequality.  However, since quaternionic K\"ahler manifolds are automatically Einstein, the same type of questions are not interesting for this class of manifolds.

{\bf Acknowledgement.} This work was done when the third author was
visiting the University of California, Irvine. He wishes to thank
the institution  for its hospitality. He also would like to thank
Professor J. Berndt for pointing out the paper of \cite{G} to him.

\section{Preliminaries on quaternionic K\"ahler manifolds}
In this section, we will recall basic properties of quaternionic K\"ahler
manifolds that will be needed in the sequel. These properties were proved by Berger \cite{B} and Ishihara \cite{I} (also see \cite{Be}).

Let $(M^n,g)$ be a Riemannian manifold, $TM$ the tangent space of $M$  and $\nabla$ the Levi-Civita
connection. The Riemannian curvature $R : TM \otimes TM \otimes TM
\longrightarrow  TM$ is defined by
\[
R(X,Y)Z = \nabla_X \nabla_Y Z - \nabla_Y \nabla_X Z - \nabla_{[X,Y]}
Z
\]
If  $\{ e_1, \cdots, e_n \}$ is an orthonormal basis of $TM$, the components of curvature tensor is defined by
\[
R_{ijkl}=\la R(e_i,e_j)e_l, e_k \ra,
 \]
the Ricci curvature is defined by
\[
\Ric_M(X,Y) = \sum_{i=1}^n \la R(X, e_i)e_i, Y\ra,
\]
and the scalar curvature is defined by
\[
{\rm S}_M = \sum_{i,j=1}^n \la R(e_i, e_j)e_j, e_i \ra.
\]

\begin{defn}
A quaternionic K\"ahler manifold $(M,g)$ is a Riemannian manifold
with a rank $3$ vector bundle $V \subset End(TM)$ satisfying
\begin{enumerate}
\item[$(a)$] In any coordinate neighborhood $U$ of $M$, there exists
a local basis $\{I,J,K\}$ of $V$ such that
 \begin{eqnarray*}
 & & I^2=J^2=K^2 = - 1 \\
 & & IJ=-JI = K \\
 & & JK = -KJ =I \\
 & & KI=-IK =J
 \end{eqnarray*}
and
\[
\la IX,IY \ra=\la JX,JY \ra=\la KX,KY \ra =\la X,Y \ra
\]
for all $X,Y \in TM$.
\item[$(b)$] If $\phi \in \Gamma(V)$, then $\nabla_X \phi \in
\Gamma(V)$ for all $X \in TM$.
\end{enumerate}
\end{defn}

\begin{rema}
It follows from $(a)$ that $\dim M = 4n$. A well known fact about
$4n$-dimensional Riemannian manifold is  that it is quaternionic K\"ahler if and only if
its restricted holonomy group is contained in $Sp(n)Sp(1)$.
\end{rema}

The $4$-dimensional Riemannian manifolds with holonomy $Sp(1)Sp(1)$
are simply the oriented Riemanian manifolds, naturally  we only
consider those when $n\ge 2$.

Notice that in general $I,J,K$ are not defined everywhere on $M$.
For example, the  canonical quaternionic projective space $QP^n$
admits no almost complex structure.

On the other hand, the vector space generated by $I,J,K$  is well
defined at each point of $M$ and this $3$-dimensional subbundle $V$
of $\End(TM)$ is in fact ``globally parallel'' under the Levi-Civita
connection $\nabla$  of $g$. A basic fact about the connection is
the following lemma.

\begin{lem}
The condition $(b)$ is equivalent to the following condition:
 \begin{align*}
 \nabla_X I &= c(X) J -b(X) K, \\
 \nabla_X J &= -c(X) I + a(X) K, \\
 \nabla_X K &= b(X) I -a(X) J,
 \end{align*}
where $a,b,c$ are local $1$-forms.
\end{lem}

\begin{defn}
Let $(M,g)$ be a quaternionic K\"ahler manifold. We can  define a
$4$-form by
\[
\Omega = \omega_1 \wedge \omega_1 +  \omega_2 \wedge \omega_2 +
\omega_3 \wedge \omega_3,
\]
where
 \begin{align*}
 \omega_1 &= \la \cdot, I\cdot \ra, \\
 \omega_2 &= \la \cdot, J\cdot \ra ,\\
 \omega_3 &= \la \cdot, K\cdot \ra.
 \end{align*}

\end{defn}

Let $\{ e_1,Ie_1,Je_1,Ke_1, \cdots, e_n,   Ie_n, Je_n, Ke_n \}$ be
an orthonormal basis of $TM$ and $\{ \theta^1,  I\theta^1,
J\theta^1, K\theta^1,\cdots, \theta^n,  I\theta^n, J\theta^n,
K\theta^n, \}$ the dual basis. It follows that
 \begin{align*}
 \omega_1 &= \sum_{i=1}^n \left( \theta^i \wedge I\theta^i + J\theta^i \wedge K \theta^i \right), \\
 \omega_2 &= \sum_{i=1}^n \left( \theta^i \wedge J\theta^i + K\theta^i \wedge I \theta^i \right), \\
 \omega_3 &= \sum_{i=1}^n \left( \theta^i \wedge K\theta^i + I\theta^i \wedge J \theta^i
 \right),
 \end{align*}
and
 \begin{align*}
 \Omega =& \sum_{i,j} \left( \theta^i \wedge I\theta^i \wedge \theta^j \wedge I\theta^j
 + \theta^i \wedge J\theta^i \wedge \theta^j \wedge J\theta^j + \theta^i \wedge K\theta^i \wedge \theta^j \wedge K\theta^j \right)  \\
  &  + \sum_{i,j} \left( J\theta^i \wedge K\theta^i \wedge J\theta^j \wedge K\theta^j
 + K \theta^i \wedge I\theta^i \wedge K \theta^j \wedge I \theta^j \right. \\
 &  +\left. I \theta^i \wedge J \theta^i \wedge I \theta^j \wedge J \theta^j \right)  \\
  &  + 2 \sum_{i,j} \left( \theta^i \wedge I\theta^i \wedge J \theta^j \wedge K \theta^j
 + \theta^i \wedge J\theta^i \wedge K \theta^j \wedge I \theta^j \right. \\
  &  + \left. \theta^i \wedge K\theta^i \wedge I \theta^j \wedge J \theta^j
 \right).
 \end{align*}

\begin{lem}\label{lemcon}
The condition $(b)$ is equivalent to the following condition:
 \begin{align*}
 \nabla_X \omega_1 &= c(X) \omega_2 -b(X) \omega_3, \\
 \nabla_X \omega_2 &= -c(X) \omega_1 + a(X) \omega_3, \\
 \nabla_X \omega_3 &= b(X) \omega_1 -a(X) \omega_2.
 \end{align*}
where $a,b,c$ are local $1$-forms.
\end{lem}

\proof{} It follows from the identities
 \begin{align*}
 (\nabla_X \omega_1)(Y,Z) &= \la Y, (\nabla_X I)Z \ra, \\
 (\nabla_X \omega_2)(Y,Z) &= \la Y,  (\nabla_X J)Z \ra, \\
 (\nabla_X \omega_3)(Y,Z) &= \la Y, (\nabla_X K)Z \ra.
 \end{align*}
 \QED

Using this lemma, we have that

\begin{thm}
The condition $(b)$ is equivalent to that $\Omega$ is parallel, that
is
\[
\nabla_X \Omega = 0
\]
for any $X \in TM$.
\end{thm}

In the following, we shall study the curvature of quaternionic
K\"ahler manifold. First we have the following lemma.

\begin{lem}\label{bracket}
If $(M^{4n},g)$ is a quaternionic K\"ahler manifold, then
 \begin{align*}
 \left[R(X,Y), I \right] &= \gamma(X,Y) J - \beta(X,Y) K, \\
 \left[R(X,Y), J \right] &= -\gamma(X,Y) I + \alpha(X,Y) K, \\
 \left[R(X,Y), K \right] &= \beta(X,Y) I -\alpha(X,Y) J,
 \end{align*}
where $\alpha$, $\beta$ and $\gamma$ are local $2$-forms given by
 \begin{align*}
 \alpha &= da + b \wedge c, \\
 \beta &= db + c \wedge a ,\\
 \gamma &= dc + a \wedge b.
 \end{align*}
\end{lem}

\begin{cor}\label{curvature}
If $(M^{4n},g)$ is a quarternionic K\"ahler manifold, then
 \begin{align*}
\la R(X,Y) Z, I Z \ra + \la R(X,Y) J Z, KZ\ra &= \alpha (X,Y)\, |Z|^2, \nonumber \\
\la R(X,Y) Z, J Z \ra + \la R(X,Y) K Z, IZ \ra&= \beta (X,Y)\, |Z|^2, \nonumber \\
\la R(X,Y) Z, K Z\ra + \la R(X,Y) I Z, JZ \ra &= \gamma (X,Y)\,
|Z|^2.
 \end{align*}
\end{cor}

The following lemma is the key for quaternionic K\"ahler manifolds.

\begin{lem}
If $(M^{4n},g)$ is a quaternionic K\"ahler manifold and $n \geq 2$,
then
 \begin{equation}\label{eqn1}
    \alpha(X,IY)=\beta(X,JY)=\gamma(X,KY)=-\frac{1}{n+2}\Ric_M (X,Y).
 \end{equation}

\end{lem}

As applications of the above lemma, one can show the following two
main theorems on curvature of quaternionic K\"ahler manifolds.

\begin{thm}
If $(M^{4n},g)$ is a quaternionic K\"ahler manifold and $n \geq 2$,
then $(M^{4n},g)$ is Einstein, that is, there is a constant $\delta$
such that
\[
\Ric_M (g) = 4(n+2)\delta g.
\]
\end{thm}

\begin{thm}\label{thmcur}
If $(M^{4n},g)$ is a quaternionic K\"ahler manifold and $n \geq 2$,
then
\begin{enumerate}
\item[$(1)$] For any tangent vector $X$, the sectional curvature
satisfies
\begin{align*}
 \la R(X,IX)IX,X \ra + \la R(X,JX)JX,X\ra &  \\
  + \la R(X,KX)KX,X\ra &= 12
\delta \,|X|^4.
\end{align*}

\item[$(2)$] For any tangent vector $Y$ satisfying
\[
\la Y,X\ra = \la Y,IX\ra = \la Y,JX\ra = \la Y, KX\ra =0,
\]
the sectional curvature satisfies
\begin{align*}
\la R(X,Y)Y,X\ra+\la R(X,IY)IY,X\ra +  &\\
\la R(X,JY)JY,X\ra+ \la R(X,KY)KY,X\ra &= 4 \delta \,|X|^2\,|Y|^2,
\end{align*}
\end{enumerate}
where $4(n+2)\delta$ is the Einstein constant.
\end{thm}

Finally, we end this section with the following lemma.
\begin{lem}\label{lemTrace} Let $\gamma:[a,b]\to M $ be a geodesic with unit speed.
 If  $S=16n(n+2)\delta$, and $X_I(t), X_J(t), X_K(t)$  are parallel vector fields along $\gamma$
 such that $X_I(a)=I\gamma'(a), X_J(a)=J\gamma'(a), X_K(a)=K\gamma'(a)$,
 then
 \[ \K(\gamma'(t), X_I(t))+\K(\gamma'(t), X_J(t))+\K(\gamma'(t), X_K(t))= 12\delta,\]
 for all $t$ and $\gamma$.

 Let $Y$ be a tangent vector at $\gamma(a)$ satisfying $\la \gamma'(a), Y \ra=0$, $\la I\gamma'(a), Y\ra = 0$, $\la J\gamma'(a), Y\ra = 0$, and $\la K\gamma'(a), Y\ra  = 0.$  If we denote the parallel vector fields $Y(t),$ $Y_I(t)$, $Y_J(t)$, and $Y_K(t)$ along $\gamma$ with initial data $Y(a)= Y$, $Y_I(a) = IY,$ $Y_J(a)= JY$, and $Y_K(a)= KY$, respectively, then
  \[ \K(\gamma'(t), Y(t))+\K(\gamma'(t), Y_I(t))+\K(\gamma'(t), Y_J(t))+\K(\gamma'(t), Y_K(t))= 4\delta,\]
  for all $t$ and $\gamma.$
 \end{lem}
 {\bf Proof.}  By the discussion above, we know the
 $3$-dimensional vector space $E(t)$ spanned by $X(t), Y(t), Z(t)$
 does not depend on the choice of $I,J,K$. Hence it is parallel under
 the Levi-Civita connection. We consider
 $\la R(\cdot,\gamma'(t))\gamma'(t),\cdot\ra$ as a symmetric bilinear
 form on $E(t)$. Then $ \K(\gamma'(t), X(t))+\K(\gamma'(t), Y(t))+\K(\gamma'(t),
 Z(t))$ is its trace on $E(t)$ which independent of the choice of
 orthonormal basis. By the computation above it is equal to
 $12\delta$.  The same argument also applies to the second part of the lemma.
 \QED

\section{Laplacian comparison theorem}
For a complete Riemannian manifold $M$ and $p\in M$, let us denote  the cut
locus with respect to $p$ by $Cut(p)$.
\begin{thm}\label{thmLaplace}Let $(M^{4n},g)$ be a complete quaternionic K\"ahler manifold with  scalar curvature ${\rm S}_M\ge 16n(n+2)\delta$ and let
$r(x)$ be the distance function to a fixed point $p\in M$. Then, for
$x\notin Cut(p)$,

\begin{equation}\label{eqnLaplace-1}
    \Delta r(x)\le \left\{
  \begin{array}{ll}
     6\coth 2r(x) +4(n-1)\coth r(x) & \hbox{ when $\delta=-1$} \\
    (4n-3)r^{-1}(x) & \hbox{ when $\delta=0$} \\
    6\cot 2r(x) +4(n-1)\cot r(x) & \hbox{ when $\delta=1$.}
  \end{array}
\right.
\end{equation}

\end{thm}
{\bf Proof.}  Let $\gamma$ be the minimizing geodesic joining $p$ to $x$.  At
$x$, we choose  $\{ e_1, e_2, \cdots, e_n\}$, and two local almost complex
structures $I,J$ and $K=IJ$ such that $e_1=\nabla r$ and
\[\{e_1,Ie_1,Je_1,Ke_1, e_2,Ie_2,Je_2,Ke_2,
\cdots, e_n,Ie_n,Je_n,Ke_n\} \]
 is an orthonormal frame.  By
parallel translating along $\gamma$ we obtain an orthonormal frame
with $e_1=\nabla r$.  For convenience sake, we denote this frame by
$\{ \varepsilon_1, \varepsilon_2, \cdots, \varepsilon_{4n}\}$. Since
$|\nabla r|^2=1$ on $M\backslash Cut(p)$, by taking covariant
derivative of this equation, we have
\begin{align}
0 &= |\nabla r|^2_{kl} \nonumber\\
  &= 2\sum_{i=1}^{4n}r_{ik}r_{il}+2\sum_{i=1}^{4n}r_ir_{ikl},
 \end{align}
for each $k,l=2,\cdots, 4n.$ Since
\begin{equation*}
    r_{ikl}=r_{kli}+\sum_{j=1}^nR_{jkil}r_j,
\end{equation*}
with $R_{ijkl}=\la R(\varepsilon_i,\varepsilon_j)\varepsilon_l, \varepsilon_k\ra,$
and $r_1=1$, $r_j=0,$  $j=2,\cdots, 4n$, we have
\begin{equation}\label{eqnriccati0}
   \sum_{i=1}^{4n}r_{ik}r_{il}+r_{kl1}+R_{1k1l}=0.
\end{equation}
In particular, if $k=l$, we have
\begin{equation}\label{eqnriccati}
   \sum_{i=1}^{4n}r_{ik}^2+r_{kk1}+\K(\varepsilon_1,\varepsilon_k)=0,
\end{equation}
where $\K(\varepsilon_1,\varepsilon_k)=R_{1k1k}$ is the sectional
curvature of the 2-plane section spanned by $\varepsilon_1, \varepsilon_k$.
Using the inequality
\begin{equation*}
    \sum_{k=2}^4r_{ik}^2\ge \frac
    13(\sum_{k=2}^4r_{kk})^2,
\end{equation*}
and setting $f(t)=\sum_{k=2}^4r_{kk}$, $(\ref{eqnriccati})$ implies
that
\begin{equation}\label{eqnricatti1}
    f'(t)+\frac
    13f^2(t)+\sum_{k=2}^4\K(\varepsilon_1,\varepsilon_k)\le 0.
\end{equation}

By Lemma \ref{lemTrace}, we have
\begin{equation}\label{eqnricatti2}
f'(t)+\frac
    13f^2(t)+12\delta\le 0.
\end{equation}
Since a smooth Riemannian metric is locally Euclidean, then
$\lim_{t\to 0}tf(t)=3$.  By a standard comparison argument for ordinary differential equations, we conclude that
\begin{equation}\label{eqnr3}
    f(t)\le\left\{
             \begin{array}{ll}
               6\cot 2t & \hbox{when $\delta=1$} \\
               3t^{-1} & \hbox{when $\delta=0$} \\
               6\coth 2t & \hbox{ when $\delta=-1$.}
             \end{array}
           \right.
\end{equation}

Similarly, using the inequality
\begin{equation*}
    \sum_{k=4i+1}^{4i+4}r_{ik}^2\ge \frac
    1{4}(\sum_{k=4i+1}^{4i+4}r_{kk})^2
\end{equation*}
for $1\le i \le n-1,$ and setting
$h_i(t)=\sum_{k=4i+1}^{4i+4}r_{kk}$, $(\ref{eqnriccati})$ implies
that
\begin{equation}\label{eqnricatti4}
    h_i'(t)+\frac
    1{4}h_i^2(t)+\sum_{k=4i+1}^{4i+4}\K(\varepsilon_1,\varepsilon_k)\le 0.
\end{equation}
Together with Lemma 1.5 asserting that
\[\sum_{k=4i+1}^{4i+4}\K(\varepsilon_1,\varepsilon_k)= 4\delta,\]
we have
\begin{equation}\label{eqnricatti5}
h_i'(t)+\frac
    1{4}h_i^2(t)+4\delta\le 0.
\end{equation}
Hence, as before, we conclude that
\begin{equation}\label{eqnr4}
    h_i(t)\le\left\{
             \begin{array}{ll}
               4\cot t & \hbox{when $\delta=1$} \\
               4t^{-1} & \hbox{when $\delta=0$} \\
               4\coth t & \hbox{ when $\delta=-1$.}
             \end{array}
           \right.
\end{equation}
The result follows from the equation $\Delta r(x)=f(r(x))+\sum_{i=1}^{n-1}h_i(r(x))$.\QED
\begin{rema}
The estimate in Theorem \ref{thmLaplace} is sharp since the right
hand sides are exactly the Laplacian of the distance functions of
quaternionic hyperbolic space $\Bbb{QH}^n$, quaternionic Euclidean space
$\Bbb{Q}^n$ and  quaternionic projective space $\Bbb{QP}^n$ respectively.
\end{rema}
\begin{rema}
We actually proved the estimate for Hessian of the distance
function.  In particular,
\begin{equation}\label{eqnhess1}
    \sum_{k=2}^4 r_{k k}\le\left\{
             \begin{array}{ll}
               6\cot 2t & \hbox{when $\delta=1$} \\
               3t^{-1} & \hbox{when $\delta=0$} \\
               6\coth 2t & \hbox{ when $\delta=-1$.}
             \end{array}
           \right.
\end{equation}
Also for $1\le i \le n-1$, we have
\begin{equation}\label{eqnhess2}
    \sum_{k=4i+1}^{4i+4} r_{kk}\le\left\{
             \begin{array}{ll}
               4\cot 2t & \hbox{when $\delta=1$} \\
               4t^{-1} & \hbox{when $\delta=0$} \\
               4\coth 2t & \hbox{ when $\delta=-1$.}
             \end{array}
           \right.
\end{equation}
\end{rema}

\begin{cor}\label{corQHA}Let $(M^{4n},g)$ be a complete quaternionic K\"ahler manifold  with  scalar curvature ${\rm S}_M\ge-16n(n+2)$. Then for any point $x\in M$
and $r>0$, the area $A(r)$ of the geodesic spheres centered at $x$ satisfies
\begin{equation}\label{eqnarea-1}
\frac{A'(r)}{A(r)}\le 6\coth 2r +4(n-1)\coth r.
\end{equation}
In particular,  $A(r)\le C(\sinh 2r)^3(\sinh r)^{4(n-1)}\le Ce^{(4n+2)r}$.

\end{cor}

\begin{cor}\label{corQHV}Let $(M^{4n},g)$ be a complete quaternionic K\"ahler manifold  with  scalar curvature ${\rm S}_M\ge-16n(n+2)$. Then for any point $x\in M$
and $0<r_1\le r_2$, the volume of the geodesic balls centered at $x$ satisfies
\begin{equation}\label{eqnvol-1}
    \frac{V_x(r_2)}{V_x(r_1)}\le
    \frac{V_{\Bbb{QH}^n}(r_2)}{V_{\Bbb{QH}^n}(r_1)},
\end{equation}
where $V_{\Bbb{QH}^n}(r)$ denotes the volume of the geodesic ball of
radius $r$ in $\Bbb{QH}^n$. In particular,  $\lambda_1(M)\le (2n+1)^2$.

\end{cor}

\begin{cor}\label{corQP}Let $(M^{4n},g)$ be a complete quaternionic K\"ahler manifold with scalar curvature ${\rm S}_M\ge 16n(n+2)$ . Then it is compact, and the diameter $d(M)\le \frac{\pi}2$,
which is the diameter of the model space $\Bbb{QP}^n$. Moreover, the
volume of $M$  is bounded by
\begin{equation}\label{eqnvol-2}
    V(M)\le V(\Bbb{QP}^n),
\end{equation}
where $V_{\Bbb{QP}^n}$ is the volume of $\Bbb{QP}^n$.
\end{cor}

\section{Quaternionic harmonicity}
In this section we will derive an over-determined system of harmonic
functions with finite  Dirichlet integral on a manifold with a
parallel form. This result was first proved by Siu \cite{S} for
harmonic maps in his proof of the rigidity theorem for K\"ahler
manifolds. Corlette \cite{C1}  gave a more systematic approach for
harmonic map with finite energy from a finite-volume quaternionic
hyperbolic space or Cayley hyperbolic plane to a manifold with
nonpositive curvature. In \cite{L}, the second author generalized
Siu's argument to harmonic functions with finite Dirichlet integral
on a K\"ahler manifold.  We will provide an argument that
generalizes Corlette's  argument to harmonic functions with finite
Dirichlet integral on a complete manifold with a parallel form. We
believe that it should be of independent interest.
\begin{thm}\label{thmrigid}Let $M$ be a complete Riemannian manifold with a parallel $p$-form $\Omega.$ Assume that $f$ is a harmonic function
with its Dirichlet integral over geodesic balls centered at $o$ of
radius $R$ satisfying the growth condition
\[
\int_{B_o(R)} |\nabla f|^2 dv = o(R^2)
\]
as $R \to \infty$,
then $f$ satisfies
\begin{equation}\label{eqnrigid1}
    d*(df\wedge\Omega)=0.
\end{equation}

\end{thm}

Before we prove the theorem, let us first recall the following operators and
some of the basic properties. For an oriented real vector space $V$ with an
inner product, we have the Hodge star operator
\begin{equation*}
    *:  \wedge^pV\to \wedge^{n-p} V.
\end{equation*}
For any $\theta\in \wedge^1V$ and $v\in V$, we also have exterior multiplication
and interior product operators
\begin{eqnarray*}
                    \varepsilon(\theta):&\wedge^pV\to\wedge^{p+1}V, \nonumber \\
                   \ell(v):&\wedge^pV\to \wedge^{p-1}V.\nonumber
                 \end{eqnarray*}
For $\theta\in \wedge^1V$ and $v\in V$ is the dual of $\theta$ by
the inner product,  if $\xi\in \wedge^pV$ we list the following
identities among the operators:
\begin{enumerate}
  \item $**\xi=(-1)^{p(n-p)}\xi,$
  \item $*\varepsilon(\theta)\xi=(-1)^p\ell(v)* \xi,$
  \item $\varepsilon(\theta)*\xi=(-1)^{p-1}*\ell(v)\xi,$
  \item $*\varepsilon(\theta)*\xi=(-1)^{(p-1)(n-p)}\ell(v)\xi,$
  \item
  $\ell(v)\varepsilon(\theta')\xi+\varepsilon(\theta)\ell(v')\xi=0,$
  where $v\bot v'$,
  \item
  $\ell(v)\varepsilon(\theta)\xi+\varepsilon(\theta)\ell(v)\xi=\xi.$
\end{enumerate}

We are now ready to prove Theorem \ref{thmrigid}.

\noindent {\bf Proof of Theorem \ref{thmrigid}.}
Let
$\eta:[0,+\infty)\to \mathbb{R} $ be a smooth function satisfying
$\eta'(t) \le 0$, and
\[ \eta(t)=\left\{
            \begin{array}{ll}
              1 & \hbox{ when $t\in [0,1]$} \\
              0 & \hbox{ when $t\in [2,+\infty]$.}
            \end{array}
          \right.
\]

For $R \ge 1$, we define the cut-off function
$\phi_R(x)=\eta(r(x)/R)$, where $r(x) $ is the distance function
from a fixed point $o\in M$,  then there is a positive constant
$C_1$ depending on $\eta$ and $C$ such that
\[
  |\nabla\phi_R(x)|\le C_1\,R^{-1}.
  \]

Since $d^2=0$, then
\begin{align}\label{eqn6.5}
  0 =& \int_Md\left\{\phi_R^2*(df\wedge\Omega)\wedge d*(df\wedge*\Omega)\right\} \nonumber\\
   =& \int_Md(\phi_R^2)\wedge*(df\wedge \Omega)\wedge d*(df\wedge*\Omega)\nonumber \\
   & +\int_M\phi_R^2\,d*(df\wedge \Omega)\wedge
   d*(df\wedge*\Omega).
\end{align}
We claim that
\begin{equation}\label{eqn7}
*d*(df\wedge \Omega)=(-1)^{n-1}d*(df\wedge *\Omega).
\end{equation}
In fact, for any  point $x\in M$, we can choose an orthonormal
tangent basis $\{e_i\}_{i=1}^m$ in a neighborhood of $x$ such that
$\nabla_{e_i}e_j(x)=0$. Denote by $\{\theta^i\}_{i=1}^m$ the dual
basis of $\{e_i\}_{i=1}^m$. Then for $\omega\in
\wedge^p(T^*M)$ we have
\begin{equation*}
    d\omega=\varepsilon(\theta^i)\nabla_{e_i} \omega.
\end{equation*}
Hence
\begin{align*}
  d*(df\wedge * \Omega) &= d*\varepsilon(df)*\Omega\nonumber \\
   &= (-1)^{(p-1)(m-p)}d[\ell(\nabla f)\Omega] \nonumber\\
   &= (-1)^{(p-1)(m-p)} \sum_{i=1}^m\varepsilon(\theta_i)\nabla_{e_i}(\ell(\nabla f)\Omega)\nonumber\\
   &=
   (-1)^{(p-1)(m-p)}\sum_{i,j=1}^m\varepsilon(\theta_i)(\nabla_{e_i}\nabla_{e_j}f)(\ell(e_j)\Omega)\nonumber\\
   &=
   (-1)^{(p-1)(m-p)}\sum_{i,j=1}^mf_{ij}\varepsilon(\theta_i)(\ell(e_j)\Omega),
\end{align*}
where $f_{ij}=\nabla_{e_i}\nabla_{e_j}f$ and the facts $\Omega$ is
parallel and $\nabla_{e_i}e_j(x)=0$ have been used. On the other
hand,

\begin{align}\label{eqn7.5}
  *d*(df\wedge  \Omega) &= *d*\varepsilon(df)\Omega\nonumber \\
   &= *\sum_{i=1}^m\varepsilon(\theta_i)\nabla_{e_i}(*\varepsilon(\sum_{j=1}^mf_j\theta_j)d[\ell(\nabla f)\Omega] \nonumber\\
   &=
   *\sum_{i,j=1}^mf_{ij}\varepsilon(\theta_i)*\varepsilon(\theta_j)(\Omega)\nonumber\\
   &=
   (-1)^{p(m-p-1)}\sum_{i,j=1}^mf_{ij}\ell(e_i)\varepsilon(\theta_j)\Omega\nonumber\\
   &=
   (-1)^{p(m-p-1)}\left(\sum_{i=1}^mf_{ii}\ell(e_i)\varepsilon(\theta_i)\Omega+\sum_{i\ne
   j}^mf_{ij}\ell(e_i)\varepsilon(\theta_j)\Omega\right)\nonumber\\
   &=
   (-1)^{p(m-p-1)}\left(\sum_{i=1}^mf_{ii}[\Omega-\varepsilon(\theta_i)\ell(e_i)\Omega]-\sum_{i\ne
   j}^mf_{ij}\varepsilon(\theta_j)\ell(e_i)\Omega\right)\nonumber\\
   &=
   (-1)^{p(m-p-1)}\sum_{i,j=1}^mf_{ij}\varepsilon(\theta_i)(\ell(e_j)\Omega),
\end{align}
where we used $f_{ij}=f_{ji}$ and $\sum_{i=1}^mf_{ii}=0$. So the
claim is proved. By $(\ref{eqn6.5})$, we have
\begin{align}\label{eqn8}
     & \int_M \phi_R^2\,|d*(df\wedge\Omega)|^2dv\nonumber \\
    & = (-1)^m\int_Md(\phi_R^2)\wedge*(df\wedge \Omega)\wedge
    d*(df\wedge*\Omega)\nonumber \\
    & \le 2\left(\int_M \phi_R^2\,|d*(df\wedge\Omega)|^2dv\right)^{\frac 12} \left( \int_M |d\phi_R|^2|*(df\wedge \Omega)|^2dv \right)^{\frac 12}
\end{align}
On the other hand, $(\ref{eqn7})$ and the fact that $\omega$ is
bounded imply that there exists a constant $C_2>0$, such that
\begin{align*}
|*(df \wedge \Omega)| &\le C_2\,|df|\\
    |d*(df\wedge*\Omega)|&= |d*(df \wedge \Omega)|.
\end{align*}
Hence combining with $(\ref{eqn8})$ and using the definition of
$\phi_R$ we conclude that
\[
\int_{B_o(R)}|d*(df\wedge\Omega)|^2dv\le C_1R^{-2} \int_{B_o(2R)}
|df|^2dv.
\]
The assumption on the growth of the Dirichlet integral of $f$ implies that the right hand side tends to zero as $R \to \infty.$  Therefore $d*(df\wedge\Omega)=0$, and the proof is complete.
\QED

\begin{lem}
Let $(M^{4n},g)$ be a quarternionic K\"ahler manifold and $n \geq
2$. If $f$ is a  function on $M$ satisfying
\begin{equation}\label{eqn9}
d* ( df \wedge \Omega) =0
\end{equation}
 for the $4$-form $\Omega$   determined by the quaternionic K\"ahler
structure, then $f$ is quaternionic harmonic, namely,  for any nonzero tangent vector $X$,
\[
f_{X,X} + f_{IX,IX} + f_{JX,JX} + f_{KX,KX} = 0
\]
where $f_{X,X} = \nabla df (X,X)$.
\end{lem}

 \proof{}
 Let
 \begin{align*}
 \{ e_A \}_{A=1}^{4n} &= \{ e_1,  e_2, \cdots,e_n,I e_1,  Ie_2, \cdots,Ie_n,\\
 &  \qquad Je_1,  Je_2, \cdots,Je_n,Ke_1,  Ke_2, \cdots,Ke_n \}
 \end{align*}
 be an orthonormal basis of $TM$ and $\{\omega_A\}$ the dual basis with $e_1=\frac{X}{\|X\|}$. Since $\Omega$ is parallel, by $(\ref{eqn7.5})$
 and $(\ref{eqn9})$, we have

 \begin{align*}
 0 &= \sum_{A=1}^{4n} (\nabla_{e_A} d
 f) \wedge \ell(e_A) \Omega \\
 &= \sum_{A, B=1}^{4n} f_{e_A, e_B} \,\omega_B \wedge \ell(e_A) \Omega
 \end{align*}
 where we have used the fact that $f$ is a harmonic function. Hence equation $(\ref{eqn9})$ implies
 \[
 \sum_{A,B=1}^{4n} f_{e_A, e_B} \, \omega_B \wedge \ell(e_A) \Omega = 0
 \]
 Comparing the coefficient of $\omega_i \wedge I\omega_i \wedge J \omega_i \wedge K\omega_i$ on both sides by
 the explicit formula for $\Omega$ given before, we obtain that
 \[
 6 \left( f_{e_i,e_i} + f_{Ie_i,Ie_i} + f_{Je_i,Je_i} + f_{Ke_i,Ke_i} \right) = 0
 \]
 for all $e_i$, $(1 \leq i \leq n)$. So the proof is complete.
 \QED

The following corollary is an immediate consequence of the lemma.
\begin{cor}\label{corQH} Let $M^{4n}$ be a complete quaternionic K\"ahler  manifold.  Assume that $f$ is a harmonic function
with its Dirichlet integral satisfying the growth condition
\[
\int_{B_o(R)} |\nabla f|^2 dv =o(R^2)
\]
as $R \to \infty$, then $f$ must satisfy
\begin{equation}\label{eqnrigid}
    d*(df\wedge\Omega)=0,
\end{equation}
where $\Omega$  is the parallel 4-form determined by the quaternionic K\"ahler
structure. Moreover, $f$ is quaternionic harmonic.
\end{cor}

 \section{Uniqueness of infinite volume end}

 Recall that for any complete manifold if $\lambda_1(M) >0$ then $M$ must be nonparabolic.  In particular, $M$ must have at least one nonparabolic ends.  It was also proved in \cite{LW1} that under the assumption that $\lambda_1(M) >0$, an end is nonparabolic if and only if it has infinite volume.

 Let us assume that $M$ has at  least two
 nonparabolic ends, $E_1$ and $E_2$. A construction of Li-Tam \cite{LT}
 asserts that one can construct a nonconstant bounded harmonic
 function with finite Dirichlet integral.  The harmonic function $f$
 can be obtained by taking a convergent subsequence of the harmonic
 functions  $f_R$, as $R\to +\infty$, satisfying
 \[ \Delta f_R=0 \qquad \textrm{ on }B(R),
 \]
 with boundary conditions
 \[
 f_R=1 \quad \hbox{on $\partial B(R)\cap E_1$}
 \]
 and
 \[
 f_R=0 \quad \hbox{on $\partial B(R)\setminus E_1$.}
\]
 It follows from the maximum principle that $0\le f_R\le 1$, hence $
 0\le f\le 1. $  We need the following estimates from \cite{LW1}(Lemma 1.1 and 1.2 in \cite{LW1}),
  and \cite{LW3}(Lemma 5.1 in \cite{LW3}).

\begin{lem}\label{lemlw1}
Let $M$ be a complete Riemannian manifold with $\lambda_1(M)>0$.
Suppose $M$ has at least two nonparabolic ends and $E$ be an end of
$M$. Then for the harmonic function $f$ constructed above, it must
satisfy the following growth conditions:
\begin{enumerate}
  \item
There exists a constant $a$ such that $f-a\in L^2 (E)$. Moreover,
the function $f-a$ must satisfy the decay estimate
\begin{equation*}
    \int_{E(R+1)\setminus E(R)}(f-a)^2\le C\exp (-2\sqrt{\lambda_1(E)}R)
\end{equation*}
for some constant $C>0$ depending on $f$, $\lambda_1(E)$ and the
dimension of $M$.
  \item The Dirichlet integral of the function $f$ must satisfy the
  decay estimate
  \begin{equation*}
    \int_{E(R+1)\setminus E(R)}|\nabla f|^2\le C\exp
    (-2\sqrt{\lambda_1(E)}R),
\end{equation*}
and
\begin{equation*}
    \int_{ E(R)}\exp (-2\sqrt{\lambda_1(E)}r(x))|\nabla f|^2\le CR
\end{equation*}
for $R$ sufficiently large.
 \end{enumerate}
\end{lem}

\begin{lem}\label{lemlw3}
Let $M$ be a complete Riemannian manifold with at least two
nonparabolic ends and $\lambda_1(M)>0$.  Then for the harmonic
function $f$ constructed above, for any $t\in (\inf f, \sup f)$ and
$(a,b)\subset (\inf f, \sup f)$,
\begin{equation*}
    \int_{\mathcal{L}(a,b)} |\nabla f|^2=(b-a)\int_{l(b)}|\nabla f|,
\end{equation*}
where
\begin{equation*}
    l(t)=\{x\in M| f(x)=t\},
\end{equation*}
and
\begin{equation*}
    \mathcal{L}(a,b)=\{x\in M| a<f(x)<b\}.
\end{equation*}
 Moreover,
\begin{equation*}
    \int_{l(t)}|\nabla f|=\int_{l(b)}|\nabla f|.
\end{equation*}
\end{lem}

We are now ready to prove Theorem \ref{thmnonp}.\\
\proof{ of Theorem \ref{thmnonp}}  Suppose to the contrary that  there exist two ends $E_1$
and $E_2$ with infinite volume. The assumption that $\lambda_1(M)>0$ implies that
they are nonparabolic. By the construction above,  there exists a
harmonic function $f$ with finite energy such that
\[
\liminf_{x\to \infty,\, x\in E_1}f(x)=1\]
and
 \[\liminf_{x\to
\infty,\, x\in E_2}f(x)=0.\]
 The Bochner formula implies that
\begin{equation}\label{eqn10}
\frac12 \Delta |\nabla f|^2=\Ric_M (\nabla f, \nabla f)+|\nabla^2
f|^2.
\end{equation}
We now choose  an orthonormal basis $\{ e_A \}_{A=1}^{4n}$  satisfying

\[
 \{ e_1,  e_2, \cdots,e_n,I e_1,  Ie_2, \cdots,Ie_n,Je_1,  Je_2, \cdots,Je_n,Ke_1,  Ke_2, \cdots,Ke_n \}
\]
 with $e_1=\frac{\nabla f}{|\nabla f|}.$
Corollary \ref{corQH} implies that
\begin{equation*}
    \sum_{i=0}^3f_{(in+1)(in+1)}=0.
\end{equation*}
Therefore, applying the arithmetic-geometric means, we have
\begin{align}
 |\nabla^2
f|^2 &= \sum_{A,B=1}^{4n}f_{AB}^2 \nonumber \\
   &\ge f_{11}^2+\sum_{i=1}^3f_{(in+1)(in+1)}^2+2\sum_{A=2}^{4n}f_{1A}^2  \nonumber\\
   &\ge f_{11}^2+\frac 13(\sum_{i=1}^3f_{(in+1)(in+1)})^2+2\sum_{A=2}^{4n}f_{1A}^2  \nonumber\\
   &\ge \frac 43|\nabla|\nabla f||^2,
\end{align}
hence combining with $(\ref{eqn10})$ we obtain
\begin{equation}\label{eqn11.5}
\frac12\Delta |\nabla f|^2\ge -\kk |\nabla f|^2+\frac
43|\nabla|\nabla f||^2.
\end{equation}
If we write $u=|\nabla f|^{\frac23}$, then
\begin{equation}\label{eqn12}
 \Delta u\ge -\frac {8(n+2)}3 u.
\end{equation}
We want to prove that the above inequality is actually an equality. The argument follows from that in \cite{LW4} after making suitable modification to fit our situation.  For
any compactly supported  smooth function  $\phi$ on $M$, we have
\begin{align}\label{eqn13}
 0&\le \int_M\phi^2u\left(\Delta u +\frac{8(n+2)}3 u \right) \nonumber \\
 &\le -2\int_M\phi u\langle    \nabla u,\nabla \phi\rangle  -\int_M\phi^2|\nabla u|^2 +\lambda_1(M)\int_M(\phi u)^2 \nonumber \\
   &\le  -2\int_M\phi u\langle    \nabla u,\nabla \phi\rangle  -\int_M\phi^2|\nabla u|^2 +\int_M|\nabla(\phi u)|^2 \nonumber \\
   &= \int_M|\nabla\phi|^2u^2.
\end{align}

Let us choose $\phi=\psi\chi$ to be the product of two compactly
supported functions. For any $\varepsilon\in
(0,\frac12)$, we define
\begin{equation*}
    \chi(x)=\left\{
             \begin{array}{ll}
               0 & \hbox{ on $\mathcal{L}(0,\sigma\varepsilon)\cup \mathcal{L}(1-\frac{\varepsilon}2,1)$} \\
               (\log 2)^{-1}(\log f-\log(\frac{\varepsilon}2)) & \hbox{on $\mathcal{L}(\frac{\varepsilon}2,\varepsilon)\cap(M\setminus E_1)$} \\
               (\log 2)^{-1}(\log (1-f)-\log(\frac{\varepsilon}2)) & \hbox{on $\mathcal{L}(1-\varepsilon,1-\frac{\varepsilon}2)\cap E_1$} \\
               1 & \hbox{otherwise.}
             \end{array}
           \right.
\end{equation*}
For $R>1$ we define
\begin{equation*}
    \psi=\left\{
           \begin{array}{ll}
             \quad 1 & \hbox{on $B(R-1)$} \\
             R-r & \hbox{on $B(R)\setminus B(R-1)$} \\
             \quad 0 & \hbox{on $M\setminus B(R)$.}
           \end{array}
         \right.
\end{equation*}
Applying to the right hand side of $(\ref{eqn13})$, we obtain
\begin{equation}\label{eqn14}
    \int_M|\nabla \phi|^2u^2\le 2\int_M |\nabla \psi|^2\chi^2|\nabla
    f|^\frac43+2\int_M |\nabla \chi|^2\psi^2|\nabla
    f|^\frac43.
\end{equation}
Since $\Ric_M \ge -4(n+2)$, then the local estimate of
Cheng-Yau \cite{CY} (see also \cite{LW2})  implies  that there exists
a constant depending on $n$ such that
\begin{equation*}
    |\nabla f|(x)\le C|1-f(x)|.
\end{equation*}
On $E_1$, the first term of $(\ref{eqn14})$ satisfies
\begin{equation}\label{eqn15}
    \int_M |\nabla \psi|^2\chi^2|\nabla f|^\frac43\le \left(\int_{\Omega}|\nabla
    f|^2\right)^{\frac 23}\left(\int_{\Omega}1\right)^{\frac 13},
\end{equation}
where $\Omega=E_1\cap (B(R)\setminus B(R-1))\cap
(\mathcal{L}(1-\varepsilon,1-\frac{\varepsilon}2)\cup\mathcal{L}(\frac{\varepsilon}2,\varepsilon)$.
Since
\begin{align*}
  \int_{\Omega}1 &\le 4\int_{\Omega}\frac{(1-f)^2}{\varepsilon^2} \\
   &\le \frac{4}{\varepsilon^2}\int_{\Omega}(1-f)^2  \\
  &\le 4C \varepsilon^{-2}\exp(-2\sqrt{\lambda_1}R),
\end{align*}
where in the last inequality we have used Lemma \ref{lemlw1}. Again
by Lemma \ref{lemlw1}, from $(\ref{eqn15})$ we have
\begin{equation}\label{eqn16}
    \int_M |\nabla \psi|^2\chi^2|\nabla f|^\frac43\le C \varepsilon^{-\frac
    23}\exp(-2\sqrt{\lambda_1}R).
\end{equation}
For the second term of $(\ref{eqn14})$ we have
 \begin{align*}
& \int_{E_1} |\nabla \chi|^2\psi^2|\nabla f|^\frac43 \\
&\le (\log 2)^{-2}\int_{\mathcal{L}(1-\varepsilon,1- \frac{\varepsilon}2)\cap E_1\cap B(R)}|\nabla f|^{\frac 43+2}(1-f)^{-2} \\
                                                &\le C(\log 2)^{-2}\int_{\mathcal{L}(1-\varepsilon,1-\frac{\varepsilon}2)\cap E_1\cap B(R)}|\nabla f|^{2}(1-f)^{-\frac23}.
                                             \end{align*}

Using the co-area formula and Lemma \ref{lemlw3} we have
\begin{align*}
 & \int_{\mathcal{L}(1-\varepsilon,1-\frac{\varepsilon}2)\cap E_1\cap B(R)}|\nabla f|^{2}(1-f)^{-\frac23} \\
 &\le \int_{1-\varepsilon}^{1-\frac{\varepsilon}2}(1-t)^{-\frac23}\int_{l(t)\cap E_1\cap B(R)}|\nabla f|dAdt \\
  &\le C\int_{l(b)} |\nabla f| dA\int_{1-\varepsilon}^{1-\frac{\varepsilon}2}(1-t)^{-\frac23}dt\\
   &= -3C[(1-t)^{\frac13}]_{1-\varepsilon}^{1-\frac{\varepsilon}2}\int_{l(b)} |\nabla f| dA\\
   &= 3C\varepsilon^{\frac13}\int_{l(b)} |\nabla f|
   dA.
\end{align*}
Combining the above inequality with $(\ref{eqn16})$ we have
\begin{equation}\label{eqn17}
    \int_{E_1}|\nabla \phi|^2u^2\le
    C(\varepsilon^{\frac23}\exp(-2\sqrt{\lambda_1}R)+\varepsilon^\frac13).
\end{equation}
A similar argument using $f$ instead of $1-f$ on the other end
yields the estimate
\[
    \int_{M\setminus E_1}|\nabla \phi|^2u^2\le
    C(\varepsilon^{\frac23}\exp(-2\sqrt{\lambda_1}R)+\varepsilon^\frac13).
\]
Letting $R\to \infty$ and $\varepsilon \to 0$, we have
\begin{eqnarray}\label{eqn19}
\Delta u= -\frac{8(n+2)}3 u
\end{eqnarray}
with $\lambda_1(M) = \frac{8(n+2)}3,$ since $f$ is nonconstant and
$u$ cannot be identically zero. Therefore all the inequalities used
to prove $(\ref{eqn12})$  are equalities. Thus there exists a
function $\mu$, such that,
\begin{equation}\label{eqn19.6}
     (f_{AB})=\left(
       \begin{array}{cccc}
         D_1 & &  &  \\
          & D_2 &  &  \\
          &  & D_2&  \\
          &  & & D_2\\
       \end{array}
     \right),
\end{equation}
where $D_1$ and $D_2$ are $n\times n$ matrices defined by
\begin{equation*}
    D_1=\left(
       \begin{array}{cccc}
         -3\mu & &  &  \\
          & 0 &  &  \\
          &  & \cdots &  \\
          &  & & 0 \\
       \end{array}
     \right)
\end{equation*}
and
\begin{equation*}
    D_2=\left(
       \begin{array}{cccc}
         \mu & &  &  \\
          & 0 &  &  \\
          &  & \cdots&  \\
          &  & & 0 \\
       \end{array}
     \right).
\end{equation*}

Since $f_{1\alpha}=0$ for $ \alpha\ne 1$ implies that $|\nabla f|$
is constant along the level set of $f$. Moreover, regularity of the
equation $(\ref{eqn19})$ implies that $|\nabla f|$  can never be
zero. Hence $M$ must be diffeomorphic to $\mathbb{R}\times N$, where
$N$ is given by the level set of $f$.  Also $N$ must be compact
since we assume that $M$ has at least 2 ends.

Fix a level set $N_0$ of $f$, consider $ (-\varepsilon, \varepsilon)
\times N_0\subset M$. Note that  $\{e_A \}$ is an orthonormal basis of
$TM$ such that $e_1$ is the normal vector to $N_0$ and $\{e_{\alpha}
\}$ are the tangent vectors of $N_0$. We shall compute the sectional
curvature
\[
\K(e_1, e_{\alpha}) = \langle  R(e_1, e_{\alpha} )e_{\alpha},
e_1\rangle  .
\]
We claim that
\[
\nabla_{e_1} e_1 = 0.
\]
Indeed it suffices to prove all integral curves $\eta(t)$ of the vector field
$e_1=\frac{\nabla f}{|\nabla f|}$ emanating from $N_0$ are geodesics. For any point
$\eta(t_0)$, let $\gamma$ be the geodesic realizing the distance
between $\eta(t_0)$ and $N_0$. Then $\gamma$ is perpendicular to
every level set $N_t$. So $\gamma'$  is parallel to $e_1$ along
$\gamma$. This implies $\gamma$ coincides with the integral curve of
$e_1$.

Let $(h_{\alpha\beta})$ with $2\le \alpha, \beta \le 4n$ be the second fundamental form of the level
set of $f$. Then
\begin{eqnarray}\label{eqn19.7}
    h_{\alpha \beta} \,f_1 = -f_{\alpha\beta},
  \end{eqnarray}
and \[ \nabla_{e_{\alpha}} e_1 = -\sum_{\beta=2}^{4n} h_{\alpha
\beta} e_{\beta}.
\]
By the definition of curvature tensor, we have
\begin{align*}
   \langle  R(e_1, e_{\alpha} )e_1, e_{\alpha}\rangle   &=
  \langle  \nabla_{e_1} \nabla_{e_{\alpha}} e_1 -  \nabla_{e_{\alpha}}\nabla_{e_1} e_1 -\nabla_{[e_1, e_{\alpha}]} e_1, e_{\alpha}\rangle   \\
   &= \langle  \nabla_{e_1} \nabla_{e_{\alpha}} e_1, e_{\alpha}\rangle   - \langle  \nabla_{[e_1,
e_{\alpha}]} e_1, e_{\alpha}\rangle   \\
   &=  \langle  \nabla_{e_1} \nabla_{e_{\alpha}} e_1, e_{\alpha}\rangle   - \langle  \nabla_{\nabla_{e_1}
e_{\alpha} - \nabla_{e_{\alpha}} e_1} e_1,
   e_{\alpha}\rangle   \\
   &= \langle  \nabla_{e_1} \nabla_{e_{\alpha}} e_1, e_{\alpha}\rangle   - \sum_{\beta=2}^{4n}\langle  \nabla_{e_1}
e_{\alpha}, e_{\beta}\rangle   \langle  \nabla_{e_{\beta}} e_1, e_{\alpha}\rangle   \\
    & ~~~~~~~~~~~~~~~~~~~~~+\sum_{\beta=2}^{4n}\langle   \nabla_{e_{\alpha}} e_1, e_{\beta}\rangle  \langle    \nabla_{e_{\beta}} e_1, e_{\alpha}\rangle   \\
   &= - \sum_{\beta=2}^{4n}\langle  \nabla_{e_1} (h_{\alpha \beta} e_{\beta}) , e_{\alpha}\rangle   +\sum_{\beta=2}^{4n} h_{\alpha \beta}
\langle  \nabla_{e_1} e_{\alpha}, e_{\beta}\rangle
+\sum_{\beta=2}^{4n} h_{\alpha \beta}^2 \\
    &= -\sum_{\beta=2}^{4n}\langle  (e_1 h_{\alpha \beta} )e_{\beta}, e_{\alpha}\rangle   -\sum_{\beta=2}^{4n} h_{\alpha \beta} \langle  \nabla_{e_1}
e_{\beta} , e_{\alpha}\rangle  \\
& ~~~~~~~~~~~~~~~~~~~+\sum_{\beta=2}^{4n} h_{\alpha \beta} \langle
\nabla_{e_1} e_{\alpha}, e_{\beta}\rangle    +
    h_{\alpha \beta}^2 \\
    &= -e_1 h_{\alpha \alpha}+2\sum_{\beta=2}^{4n}h_{\alpha \beta} \langle  \nabla_{e_1} e_{\alpha}, e_{\beta}\rangle   + h_{\alpha
\beta}^2.
\end{align*}
Therefore
\begin{equation} \label{eqn19.75}
\K(e_1, e_{\alpha}) =  e_1 h_{\alpha \alpha} -
2\sum_{\beta=2}^{4n}h_{\alpha \beta} \langle  \nabla_{e_1}
e_{\alpha}, e_{\beta}\rangle   - \sum_{\beta=2}^{4n}h_{\alpha
\beta}^2.
\end{equation}
Since  $h_{\alpha \beta}$ is diagonal, this implies that
\[
\K(e_1, e_{\alpha}) =  e_1 h_{\alpha \alpha}  - h_{\alpha
\alpha}^2.
\]
Combining with $(\ref{eqn19.6})$ and $(\ref{eqn19.7})$, we conclude
that
\begin{equation*}
    \K(e_1, e_2)=\K(e_1,I e_2)=\K(e_1,J e_2)=\K(e_1, Ke_2)=0
\end{equation*}
which implies $M$ is Ricci flat by  Theorem \ref{thmcur}. This
contradicts to the assumption that  $\lambda_1>\frac{8(n+2)}3>0$.
Therefore $M$ must have only one end with infinite volume. \QED

\section{Maximal first eigenvalue }

In this section, we will consider the case when $\lambda_1(M)$ is of
maximal value.

\proof{ of Theorem \ref{thmend}}According to Theorem \ref{thmnonp},
we know that $M$ has exactly one nonparabolic end. Suppose that $M$
has more than one end. Then  there must exist at least
 an end with finite volume. We divide the rest of the proof into several parts.  The first part
follows exactly as that in the proof of the corresponding theorem in the
K\"ahler case (Theorem 3.1) in \cite{LW5}.
For completeness sake, we will give a quick outline of it.

{\bf Part 1.}  Assume that $E_1$ is such an end with
finite volume given by $M\setminus B_p(1)$. Then we can choose a ray
$\eta: [0,+\infty) $ such that $\eta(0)=p$ and
$\eta[1,+\infty)\subset E_1$. The Busemann function corresponding to
$\gamma$ is defined by
\[ \beta(x)=\lim_{t\to +\infty}[t-d(x,\eta(t))].
\]
The Laplacian comparison theorem, Theorem \ref{thmLaplace}, asserts that
\[\Delta \beta\ge -2(2n+1)
\]
in the sense of distribution.
We define the function $f=\exp((2n+1)\beta)$, and using the fact that $|\nabla
\beta|= 1$ almost everywhere,  we have

\begin{align*}
\Delta f &= (2n+1)\exp((2n+1)\beta)\Delta \beta+(2n+1)^2\\
&\ge -(2n+1)^2f.
\end{align*}
Similar to the proof of above theorem,  we conclude that for any
compactly supported function $\phi$,
\begin{align*}
0 &\le
\int_M(\Delta f+(2n+1)^2f)f\phi^2\\
&\le \int_Mf^2|\nabla \phi|^2.
\end{align*}

By choosing the function $\phi$ to be

\begin{equation*}
    \phi=\left\{
           \begin{array}{ll}
             \quad 1, & \hbox{on $B_p(R)$;} \\
             \frac{2R-r(x)}{R}, & \hbox{on $B_p(2R)\setminus B_p(R)$}; \\
             \quad 0, & \hbox{on $M\setminus B_p(2R)$;}
           \end{array}
         \right.
\end{equation*}
we obtain
\begin{align*}
& \int_{M\cap E_1}f^2|\nabla \phi|^2 \\
&\le \frac{1}{R^2}\int_{(B_p(2R)\setminus B_p(R))\cap E_1}f^2\\
&\le \frac{1}{R^2}\sum_{i=1}^{[R]}\int_{(B_p(R+i)\setminus B_p(R+i-1))\cap E_1}f^2\\
&\le \frac{C}{R^2}\sum_{i=1}^{[R]} (V_{E_1}(R+i)-V_{E_1}(R+i-1))
\exp(2(2n+1)(R+i))
\end{align*}
where $V_{E_1}(R+i)$ denotes the volume of the set $E_1 \cap
B_p(R+i).$ On the other hand, the volume estimate in Theorem 1.4 of
\cite{LW1} implies that
\[
V_{E_1}(\infty)-V_{E_1}(R)\le C\exp (-2(2n+1)R),
\]
hence
\begin{align*}
&  V_{E_1}(R+i)-V_{E_1}(R+i-1)\\
&= V_{E_1}(\infty)-V_{E_1}(R+i-1) -(V_{E_1}(\infty)-V_{E_1}(R+i))\\
&\le C\exp(-2(2n+1)(R+i)).
\end{align*}
Therefore, we conclude that
\[ \int_{M\cap E_1}f^2|\nabla \phi|^2\le \frac {C}R.
\]

Let us now denote $E_2 = M \setminus (B_p(1) \cup E_1)$ to be the other end of $M$. When $x\in E_2$, following the argument in Theorem 3.1 of \cite{LW4}, we have
\[ \beta(x)\le -d(p,x)+2.
\]
Therefore
\begin{align*}
\int_{E_2}f^2|\nabla \phi|^2 &\le \frac{1}{R^2}\int_{(B_p(2R)\setminus B_p(R))\cap E_2}f^2\\
&= \frac{C}{R^2}\int_{(B_p(2R)\setminus B_p(R))\cap E_2}\exp(-2(2n+1)(r-2))\\
&\le \frac{C}{R}.
\end{align*}

Letting $R\to +\infty$, we conclude that
\begin{eqnarray} \label{eqn20}
 \Delta f+(2n+1)^2f=0,
\end{eqnarray}
and all inequalities used are indeed equalities and $f$ is smooth by
regularity of the equation $(\ref{eqn20})$. Moreover, $|\nabla
\beta|=1$, and
\[\Delta \beta =-2(2n+1).
\]
This implies that $M$ must be diffeomorphic to $\mathbb{R}\times N$, where $N$ is
given by the level set of $\beta.$ We choose  an orthonormal
basis $\{ e_i \}_{i=1}^{4n}$ as follows

\[
 \{ e_1,  e_2, \cdots,e_n,I e_1,  Ie_2, \cdots,Ie_n,Je_1,  Je_2, \cdots,Je_n,Ke_1,  Ke_2, \cdots,Ke_n \}
\]
 with $e_1=\nabla \beta.$
 Applying the Bochner formula to
$\beta $, we get

\begin{align*}
0&= \frac12\Delta |\nabla \beta|^2\\
&= \sum_{i,j=1}^{4n}\beta_{ij}^2+\Ric_M(\nabla \beta, \nabla \beta)+\sum_{i=1}^{4n}\beta_i(\Delta \beta)_i\\
&= \sum_{i,j=1}^{4n}\beta_{ij}^2-4(n+2).
\end{align*}
By the comparison theorem, we have,
\begin{eqnarray*}
\sum_{i=0}^3\beta_{(in+1)(in+1)}=-6.
\end{eqnarray*}
Hence
\begin{equation*}
     (\beta_{\alpha \beta})=\left(
       \begin{array}{cccc}
         D_1 & &  &  \\
          & D_2 &  &  \\
          &  & D_2&  \\
          &  & & D_2\\
       \end{array}
     \right),
\end{equation*}
where $D_1$ and $D_2$ are $n\times n$ matrices defined by
\begin{equation*}
    D_1=\left(
       \begin{array}{cccc}
         0 & &  &  \\
          & -1 &  &  \\
          &  & \cdots &  \\
          &  & & -1 \\
       \end{array}
     \right)
\end{equation*}
and
\begin{equation*}
    D_2=\left(
       \begin{array}{cccc}
         -2 & &  &  \\
          & -1 &  &  \\
          &  & \cdots&  \\
          &  & & -1 \\
       \end{array}
     \right)
\end{equation*}

 {\bf Part 2.} For a fix level set $N_0$ of $\beta $,  we consider $ (-\varepsilon,
\varepsilon) \times N_0\subset M$. Note that  $\{e_i \}$ is an orthonormal
basis of $TM$ such that $e_1$ is the normal vector to $N_0$ and
$\{e_{\alpha} \}$, for $2 \le \alpha \le 4n$, are the tangent vectors of $N_0$. We shall compute
the sectional curvature
\[
\K(e_1, e_{\alpha}) = \langle  R(e_1, e_{\alpha} )e_{\alpha},
e_1\rangle  .
\]
Since $ \nabla_{e_1} e_1 = 0$ implies that the integral curves of
$e_1$ are geodesics.
Let $(h_{\alpha\gamma})$ be the second fundamental form of the level
set of $\nabla \beta$. Then

\begin{align*}
    h_{\alpha \gamma} &= \langle  \nabla_{e_{\alpha}} e_{\gamma}, e_1\rangle   \\
     &= \langle  \nabla_{e_{\alpha}} e_{\gamma},  \nabla \beta\rangle   \\
     &= - \beta_{\alpha\gamma}
  \end{align*}
and
\begin{equation}\label{eqn20.1}
    \nabla_{e_{\alpha}} e_1 = -\sum_{\gamma=2}^{4n} h_{\alpha
\gamma} e_{\gamma}.
\end{equation}

By $(\ref{eqn19.75})$ in the proof of Theorem \ref{thmnonp}  we have
\begin{eqnarray*}
  \langle  R(e_1, e_{\alpha} )e_1, e_{\alpha}\rangle   =
  -e_1 h_{\alpha \alpha}+2\sum_{\gamma=2}^{4n}h_{\alpha \gamma} \langle  \nabla_{e_1} e_{\gamma}, e_{\beta}\rangle   +\sum_{\gamma=2}^{4n} h_{\alpha
\gamma}^2.
\end{eqnarray*}

Since $(h_{\alpha \gamma})$ are constant and  diagonal, then
\[
\K(e_1, e_{\alpha}) =   - h_{\alpha \alpha}^2.
\]

In particular, we have
\begin{equation*}
\K(e_1, e_{\alpha})=\left\{
  \begin{array}{ll}
    -4 & \hbox{ when $\alpha=in+1$, $i=1,2,3$} \\
    -1 & \hbox{ otherwise.}
  \end{array}
\right.
\end{equation*}
On the other hand, we also have
\begin{align*}
  \K(e_{n+1}, e_{2n+1})+\K(e_{n+1}, e_{3n+1}) &=-12-\K(e_{1}, e_{n+1})=-8 \\
  \K(e_{n+1}, e_{2n+1})+\K(e_{3n+1}, e_{2n+1}) &= -8\\
  \K(e_{3n+1}, e_{2n+1})+\K(e_{n+1}, e_{3n+1}) &=-8,
\end{align*}
hence
\[
\K(e_{n+1}, e_{2n+1})=\K(e_{n+1}, e_{3n+1})=\K(e_{2n+1},
e_{3n+1})=-4.
\]
Since for $\alpha=2,3,\cdots,n$,
\begin{align*}
  \K(Ie_1,e_{\alpha}) &= -\langle    R(Ie_1,e_{\alpha})Ie_1,e_{\alpha}\rangle    \\
   &= -\langle   IR(Ie_1,e_{\alpha})Ie_1,Ie_{\alpha}\rangle   \\
   &= \langle   R(Ie_1,e_{\alpha})e_1,Ie_{\alpha}\rangle  \\
   &= \langle     R(e_1,Ie_{\alpha})Ie_1,e_{\alpha}\rangle  \\
 &= \K(e_1,Ie_{\alpha}) \\
 &= -1,
\end{align*}
and $\K(Je_1,e_{\alpha})=\K(Ke_1,e_{\alpha})=-1$, we have
\[\K(e_{in+1},e_{\alpha}) =-1,
\]
for all $i=0,1,2,3$ and $\alpha\ne 1, n+1, 2n+1, 3n+1$.

Let $\K^N(e_{\alpha},e_{\gamma})$ denote the sectional curvature of the level set with induced metric. By Gaussian equation,
\[ \K^N(e_{\alpha},e_{\gamma})-\K(e_{\alpha},e_{\gamma})=h_{\alpha\alpha}h_{\gamma\gamma},
\]
it is straightforward to obtain
\[
\K^N(e_{n+1}, e_{2n+1})=\K^N(e_{n+1}, e_{3n+1})=\K^N(e_{2n+1},
e_{3n+1})=0,
\]
and
\begin{equation}\label{eqn19.8}
    \K^N(e_{in+1},e_{\alpha}) =1,
\end{equation}
for all $i=1,2,3$ and $\alpha\ne 1, n+1, 2n+1, 3n+1$.

 {\bf Part 3.} There is a
natural map $\varphi_t$ between the level sets $N_0$ and $N_t$ given
by the gradient flow of $\beta$. Since the integral curves are
geodesics,  $d\varphi_t(X)$ are Jacobi fields along corresponding
curves. Let $(N,g_0)=N_0$  with the induced metric. We can consider
$\varphi$ as a flow on $N$.  We claim that
\[
d\varphi_t|_{V_1}=e^{2t} \, \text{id}
\]
and
\[
d\varphi_t|_{V_2}=e^{t}\, \text{id},
\]
where $TN=V_1\oplus V_2$, $ V_1=\textrm{span} \{Ie_1, Je_1, Ke_1\}$
and $V_2=V_1^{\perp}$. Indeed for any point $q\in N_0$, denote
$e_1(t)=\nabla \beta(\varphi(t))$ and
$\{\varepsilon_{\alpha}(t)\}_{\alpha=2}^{4n}$ to be the parallel
transport of the orthonormal base $\{e_{\alpha}\}_{\alpha=2}^{4n}$
of $N_0$ at $q$ along $\varphi_t(q)$. Since both $V_1$ and $V_2$ are
$\varphi$-invariant, we have, in particular,
\begin{equation}\label{eqn20.9}
    \langle    \nabla_{e_1(t)}\varepsilon_{\alpha},\varepsilon_{\gamma}
    \rangle=0,
\end{equation}
when $\alpha\in \{n+1, 2n+1, 3n+1\}$, and $\gamma\notin \{n+1, 2n+1,
3n+1\}.$

Now we can compute $R_{1\alpha 1\gamma}$. Then
\begin{align}
  &  \langle  R(e_1, \varepsilon_{\alpha} )e_1, \varepsilon_{\gamma}\rangle \nonumber \\
  &= \langle  \nabla_{e_1} \nabla_{\varepsilon_{\alpha}} e_1 - \nabla_{\varepsilon_{\alpha}}\nabla_{e_1} e_1
-\nabla_{[e_1, \varepsilon_{\alpha}]} e_1, \varepsilon_{\gamma}\rangle  \nonumber  \\
   &= \langle  \nabla_{e_1} \nabla_{\varepsilon_{\gamma}} e_1, \varepsilon_{\alpha}\rangle   - \langle  \nabla_{[e_1,
\varepsilon_{\alpha}]} e_1, \varepsilon_{\gamma}\rangle  \nonumber \\
   &=  \langle  \nabla_{e_1} \nabla_{\varepsilon_{\alpha}} e_1, \varepsilon_{\gamma}\rangle   - \langle  \nabla_{\nabla_{e_1}
\varepsilon_{\alpha} - \nabla_{\varepsilon_{\alpha}} e_1} e_1, \varepsilon_{\gamma}\rangle  \nonumber  \\
   &=  \langle  \nabla_{e_1} \nabla_{\varepsilon_{\alpha}} e_1, \varepsilon_{\gamma}\rangle   - \sum_{\tau=2}^{4n}\langle  \nabla_{e_1}
\varepsilon_{\alpha}, \varepsilon_{\tau}\rangle   \langle  \nabla_{\varepsilon_{\tau}} e_1, \varepsilon_{\gamma}\rangle  \nonumber  \\
& ~~~ +\sum_{\tau=2}^{4n}\langle   \nabla_{\varepsilon_{\alpha}}
e_1,
\varepsilon_{\tau}\rangle  \langle \nabla_{\varepsilon_{\tau}} e_1, \varepsilon_{\gamma}\rangle \nonumber  \\
   &= - \sum_{\tau=2}^{4n}\langle  \nabla_{e_1} (h_{\alpha \tau} \varepsilon_{\tau}) , \varepsilon_{\gamma}\rangle   +\sum_{\tau=2}^{4n} h_{\gamma \tau}
\langle  \nabla_{e_1} \varepsilon_{\alpha}, \varepsilon_{\tau}\rangle
+\sum_{\tau=2}^{4n}
    h_{\alpha \tau}h_{\tau\gamma} \nonumber \\
    &= -e_1 h_{\alpha \gamma}    -\sum_{\tau=2}^{4n} h_{\alpha \tau} \langle  \nabla_{e_1}
\varepsilon_{\tau} , \varepsilon_{\gamma}\rangle \nonumber \\ &
~~~+\sum_{\tau=2}^{4n} h_{\gamma \tau} \langle  \nabla_{e_1}
\varepsilon_{\alpha}, \varepsilon_{\tau}\rangle +
\sum_{\tau=2}^{4n}h_{\alpha \tau}h_{\tau\gamma}.
\end{align}

We see that $(h_{\alpha\gamma})$ is diagonal and
\[h_{\alpha\alpha}=\left\{
                     \begin{array}{ll}
                       2, & \hbox{ when $\alpha=n+1,2n+1,3n+1$;} \\
                       1, & \hbox{ otherwise.}
                     \end{array}
                   \right.
\]
Therefore, when $\alpha\neq \gamma$,
\begin{align*} R_{1\alpha
1\gamma}&=-h_{\alpha\alpha}\langle\nabla_{e_1}\varepsilon_{\alpha},
\varepsilon_{\gamma}\rangle+h_{\gamma\gamma}\langle\nabla_{e_1}\varepsilon_{\alpha},
\varepsilon_{\gamma}\rangle\\
&=(h_{\gamma\gamma}-h_{\alpha\alpha})\langle\nabla_{e_1}\varepsilon_{\alpha},
\varepsilon_{\gamma}\rangle.
\end{align*}
Since $h_{\alpha\alpha}=h_{\gamma\gamma}$ when $\alpha, \gamma\in
\{n+1, 2n+1, 3n+1\}$ and $\alpha, \gamma\notin \{n+1, 2n+1, 3n+1\}$,
using $(\ref{eqn20.9})$, we have
\begin{equation*}
    R_{1\alpha 1\gamma}=0, \textrm{ for all }\alpha\ne \gamma.
\end{equation*}
Define
\begin{equation*}
    J_{\alpha}(t)=\left\{
                    \begin{array}{ll}
                      e^{-2t}\varepsilon_{\alpha}, & \hbox{ when $\alpha\in
\{n+1, 2n+1, 3n+1\}$;} \\
                      e^{-t}\varepsilon_{\alpha}, & \hbox{when $\alpha\notin
\{n+1, 2n+1, 3n+1\}$.}
                    \end{array}
                  \right.
\end{equation*}
Since
\[ \nabla_{\frac{\partial}{\partial
t}}d\varphi_t(e_{\alpha})|_{t=0}=[e_1,e_{\alpha}]=-\nabla_{e_{\alpha}}e_1,
\]
then we see that $J_{\alpha}$ satisfies the Jacobi equation and
initial conditions $J_{\alpha}(0)=e_{\alpha}$ and
$J'_{\alpha}(0)=e_{\alpha}=\nabla_{\frac{\partial}{\partial
t}}d\varphi_t(e_{\alpha})|_{t=0}.$ By the uniqueness theorem for the
Jacobi equations, we have $d\varphi_t(e_{\alpha})=J_{\alpha}.$ The
claim is proved.

{\bf Part 4.} We have now a family of metrics on $N$ written as

\begin{eqnarray*}
   ds_t^2 = e^{4t}\sum_{i=1}^3\omega_{in+1}^2+e^{2t}\sum_{i=0}^3\sum_{\alpha=2}^n\omega_{in+\alpha}^2,
 \end{eqnarray*}
and the metric of $M$ can rewritten as

\begin{equation}\label{eqn21}
ds^2 = dt^2 + e^{4t}\,\sum_{p=2}^4 \omega_p^2
+e^{2t}\,\sum_{\alpha=5}^{4n}  \omega_\alpha^2
\end{equation}
where $\{\omega_2, \omega_3, \omega_4, \dots, \omega_{4n}\}$ is the
dual coframe to $\{e_2, e_3, e_4, \dots, e_{4n}\}$ at $N_0.$ We also
choose that  $Ie_{4s-3}=e_{4s-2}$, $Je_{4s-3} = e_{4s-1}$, and
$Ke_{4s-3}=e_{4s}$ for $s=1, \dots, n,$ with $e_1= \frac{\p}{\p t}.$
In particular, the second fundamental form on $N_t$ must be a
diagonal matrix when written in terms of the basis
$\{e_i\}_{i=2}^{4n}$ with eigenvalues given by

\begin{equation}\label{eqn22}
(\la \n_{e_i} e_j, e_1\ra) =
\begin{pmatrix} 2I_3 & 0 \\0&I_{4(n-1)}
\end{pmatrix},
\end{equation}
where $I_k$ denotes the $k\times k$ identity matrix. Also, the
sectional curvatures of the sections containing $e_1$ are given by
$$
\K(e_1, e_p) = -4 \qquad \text{for} \qquad 2\le p \le 4
$$
and
$$
\K(e_1, e_\alpha) = -1 \qquad \text{for} \qquad 5\le \alpha \le 4n.
$$
The Guass curvature equation also asserts that
$$
R_{ijkl} = \bar R_{ijkl} + h_{li}h_{kj} - h_{ki}h_{lj},
$$
where $\bar R_{ijkl}$ is the curvature tensor on $N_t.$ In
particular,
\begin{equation}\label{eqn23}
R_{ijkl} = \left\{ \aligned  \bar R_{ijkl}& + \delta_{li} \delta_{kj} - \delta_{ki} \delta_{lj} \qquad \text{if} \qquad 5\le i, j, k, l \le 4n\\
\bar R_{ijkl}& + 4\delta_{li} \delta_{kj} - 4\delta_{ki} \delta_{lj}  \qquad \text{if} \qquad 2 \le i,j,k,l \le 4\\
 \bar R_{ijkl}& + 2 \qquad \text{if} \qquad 2\le i= l \le 4 \quad \text{and} \quad 5 \le k=j \le 4n\\
  \bar R_{ijkl}& + 2 \qquad \text{if} \qquad 2\le k= j \le 4 \quad \text{and} \quad5 \le i=l \le 4n\\
  \bar R_{ijkl}& -2 \qquad \text{if} \qquad 2\le i= k \le 4 \quad \text{and} \quad 5\le j=l \le 4n\\
 \bar R_{ijkl}& -2 \qquad \text{if} \qquad 2\le j= l \le 4 \quad \text{and} \quad 5\le i=k \le 4n\\
  \bar R_{ijkl}& -2 \qquad \text{if} \qquad 2\le k= i \le 4 \quad \text{and} \quad 5\le j=l \le 4n\\
 \bar R_{ijkl}& \qquad \text{otherwise.}\endaligned \right.
\end{equation}

We will now use $(\ref{eqn21})$ to  compute the curvature tensor of
$M$ and hence $N_0.$   Using the orthonormal coframe
$$
\eta_1 = \omega_1 = dt,
$$
$$
\eta_p = e^{2t} \omega_p
$$
$$
\eta_\alpha = e^t \omega_\alpha
$$
for $2 \le p \le 4$ and $5 \le \alpha \le 4n$, we obtain the first
structural equations
\begin{equation}\label{eqn24}
d\eta_1 = 0,
\end{equation}

\begin{align}\label{eqn25}
          d \eta_p &= 2e^{2t}\,\omega_1\wedge \omega_p +  e^{2t}\, \sum_{q=2}^4 \omega_{pq} \wedge \omega_q + e^{2t}\,\sum_{\alpha=5}^{4n}\omega_{p\alpha} \wedge \omega_\alpha\nonumber\\
&= -2 \eta_p \wedge \eta_1 +  \sum_{q=2}^4 \omega_{pq} \wedge
\eta_q+ e^t\,\sum_{\alpha=5}^{4n} \omega_{p\alpha} \wedge
\eta_\alpha,
         \end{align}

and
\begin{align}\label{eqn26}
d \eta_\alpha &= e^t\, \omega_1 \wedge \omega_\alpha + e^t \,\sum_{p=2}^4 \omega_{\alpha p}\wedge \omega_p + e^t\, \sum_{\beta=5}^{4n}\omega_{\alpha \beta} \wedge \omega_\beta\nonumber\\
&= -\eta_\alpha \wedge \eta_1 + e^{-t}\,\sum_{p=2}^4 \omega_{\alpha
p} \wedge \eta_p + \sum_{\beta=5}^{4n}\omega_{\alpha \beta} \wedge
\eta_\beta,
 \end{align}
where $\omega_{ij}$ are the connection forms of $N_0.$  In the above
and all subsequent computations, we will adopt the convention that
$5 \le \alpha, \beta \le 4n,$ $2 \le i, j \le 4n,$ $2\le o, p, q, r
\le 4,$  $ 2\le s, t \le n,$ and $1\le A, B \le 4n.$

Note that using the endomorphism $I$ and the fact that $\nabla I =
cJ - bK$,  we have
\begin{align*}
\omega_{ij}(X) &= \la \bar \n_{X} e_j, e_i \ra\\
&= \la I\n_{X} e_j, Ie_i\ra\\
&= \la \n_{X} Ie_j, Ie_i\ra +c(X)\,\la Je_j, Ie_i\ra - b(X) \, \la ke_j, Ie_i\ra\\
&= \la \n_{X} Ie_j, Ie_j \ra +c(X)\, \la e_j, Ke_i\ra + b(X) \la
e_j, Je_i\ra
 \end{align*}
for any tangent vector $X$ to $N_0,$ where $\bar \n$ denotes the
connection on $N_0.$  Hence we conclude that

\begin{equation}\label{eqn27}
\omega_{ij} = \omega_{I_iI_j} +c\,\la e_j, Ke_i\ra + b\, \la e_j,
Je_i\ra,
\end{equation}
 where $I_i$ denotes the index corresponding to $Ie_i =
e_{I_i}.$ Similarly, we have
$$
\omega_{ij} = \omega_{J_i J_j} + c\,\la e_j, Ke_i\ra + a\, \la e_j ,
Ie_i\ra,
$$
and
$$
\omega_{ij} = \omega_{K_i K_j} +b\,\la e_j, Je_i\ra + a\, \la e_j,
Ie_i\ra.
$$
Together with $(\ref{eqn22})$, we conclude that
$$
\omega_{2(4s-1)}(e_{4s}) = -1= -\omega_{2(4s)}(e_{4s-1}),
$$
$$
\omega_{2(4s-3)}(e_{4s-2}) = -1= -\omega_{2(4s-2)}(e_{4s-3}),
$$
for all $2\le s\le n,$ and
$$
\omega_{2\alpha}(e_\beta) = 0 \qquad \text{otherwise.}
$$
Similarly,
\begin{align*}
\omega_{2 \alpha}(e_p) &= \la \n_{e_p} e_\alpha, e_2 \ra\\
&= -\la \n_{e_p} Ie_\alpha, e_1 \ra\\
&=0.
\end{align*}

These identities imply that
\begin{align}\label{eqn28}
  \omega_{2(4s-3)} &= -\omega_{4s-2},\nonumber\\
\omega_{2(4s-2)} &= \omega_{4s-3},\nonumber\\
\omega_{2(4s-1)} &= -\omega_{4s},\nonumber\\
\omega_{2(4s)} &= \omega_{4s-1}.
\end{align}

A similar calculation using the endomorphisms $J$ and $K$ yield
\begin{align}\label{eqn29}
\omega_{3(4s-3)} &= -\omega_{4s-1},\nonumber\\
\omega_{3(4s-2)} &= \omega_{4s},\nonumber\\
\omega_{3(4s-1)} &= \omega_{4s-3},\nonumber\\
\omega_{3(4s)} &= -\omega_{4s-2},
\end{align}
and
\begin{align}\label{eqn30}
\omega_{4(4s-3)} &= -\omega_{4s},\nonumber\\
\omega_{4(4s-2)} &= -\omega_{4s-1},\nonumber\\
\omega_{4(4s-1)} &= \omega_{4s-2},\nonumber\\
\omega_{4(4s)} &= \omega_{4s-3}.
\end{align}

We claim that the connection forms are given by
\begin{align}\label{eqn31}
\eta_{1p} &=-\eta_{p1}\nonumber\\
&= 2 \eta_p \qquad \text{for} \qquad 2\le p \le 4,
\end{align}

\begin{align}\label{eqn32}
\eta_{1\alpha} &=-\eta_{\alpha 1}\nonumber\\
&= \eta_\alpha \qquad \text{for} \qquad 5 \le \alpha \le 4n,
\end{align}

\begin{equation}\label{eqn33}
\eta_{pq} = - \eta_{qp} = \omega_{pq},
\end{equation}

\begin{eqnarray}\label{eqn34}
\eta_{p \alpha}=-\eta_{\alpha p}=e^t\, \omega_{p \alpha},
\end{eqnarray}

\begin{align}\label{eqn35}
\eta_{(4s) \beta} &= -\eta_{\beta (4s)}\nonumber\\
&= \left\{ \aligned \omega_{(4s) \beta} -(1-&e^{-2t})\,\eta_2 \qquad \text{if} \qquad \beta= 4s-1\\
\omega_{(4s) \beta} +(1-&e^{-2t})\,\eta_3 \qquad \text{if} \qquad \beta= 4s-2\\
\omega_{(4s) \beta} -(1-&e^{-2t})\,\eta_4 \qquad \text{if} \qquad \beta= 4s-3\\
\omega_{(4s) \beta}& \qquad \text{if} \qquad \beta \neq 4s-1, 4s-2,
\text{or } 4s-3, \endaligned \right.
\end{align}

\begin{align}\label{eqn36}
\eta_{(4s-1) \beta} &= -\eta_{\beta (4s-1)}\nonumber\\
&= \left\{ \aligned \omega_{(4s-1) \beta} +(1-&e^{-2t})\,\eta_2 \qquad \text{if} \qquad \beta= 4s\\
\omega_{(4s-1) \beta} -(1-&e^{-2t})\,\eta_4 \qquad \text{if} \qquad \beta= 4s-2\\
\omega_{(4s-1) \beta} -(1-&e^{-2t})\,\eta_3 \qquad \text{if} \qquad \beta= 4s-3\\
\omega_{(4s-1) \beta}& \qquad \text{if} \qquad \beta \neq 4s, 4s-2,
\text{or } 4s-3, \endaligned \right.
\end{align}

\begin{align}\label{eqn37}
\eta_{(4s-2) \beta} &= -\eta_{\beta (4s-2)}\nonumber\\
&= \left\{ \aligned \omega_{(4s-2) \beta} -(1-&e^{-2t})\,\eta_3 \qquad \text{if} \qquad \beta= 4s\\
\omega_{(4s-2) \beta} +(1-&e^{-2t})\,\eta_4 \qquad \text{if} \qquad \beta= 4s-1\\
\omega_{(4s-2) \beta} -(1-&e^{-2t})\,\eta_2 \qquad \text{if} \qquad \beta= 4s-3\\
\omega_{(4s-2) \beta}& \qquad \text{if} \qquad \beta \neq 4s, 4s-1,
\text{or } 4s-3, \endaligned \right.
\end{align}

\begin{align}\label{eqn38}
\eta_{(4s-3) \beta} &= -\eta_{\beta (4s-3)}\nonumber\\
&= \left\{ \aligned \omega_{(4s-3) \beta} +(1-&e^{-2t})\,\eta_4 \qquad \text{if} \qquad \beta= 4s\\
\omega_{(4s-3) \beta} +(1-&e^{-2t})\,\eta_3 \qquad \text{if} \qquad \beta= 4s-1\\
\omega_{(4s-3) \beta} +(1-&e^{-2t})\,\eta_2 \qquad \text{if} \qquad \beta= 4s-2\\
\omega_{(4s-3) \beta}& \qquad \text{if} \qquad \beta \neq 4s, 4s-1,
\text{or } 4s-2. \endaligned \right.
\end{align}

Indeed, if we substitute $(\ref{eqn31}- \ref{eqn38})$ into the first
structural equations
$$
d \eta_A =\eta_{A1}\wedge \eta_1 + \sum_{q=2}^4\eta_{Aq} \wedge
\eta_q + \sum_{\beta=5}^{4n} \eta_{A\beta} \wedge \eta_{\beta}
$$
we obtain $(\ref{eqn24})$, $(\ref{eqn25})$, and $(\ref{eqn26})$.

To compute the curvature, we consider the second structural
equations.  In particular,
\begin{eqnarray*}
&&d \eta_{1p} - \eta_{1q} \wedge \eta_{qp} - \eta_{1 \alpha} \wedge \eta_{\alpha p} \\
&&\qquad= 2d\eta_p -2 \eta_q \wedge \eta_{qp} - \eta_{\alpha} \wedge   \eta_{\alpha p}  \\
& &\qquad= -4 \eta_p \wedge \eta_1  +\eta_\alpha \wedge \eta_{\alpha p} \\
&&\qquad= -4 \eta_p \wedge \eta_1 + e^t\, \omega_{p \alpha} \wedge
\eta_\alpha.
\end{eqnarray*}
Hence using $(\ref{eqn28} -\ref{eqn30})$, we have
\begin{align*}
    R_{1p1p} &= -4,\\
R_{12(4s-1) (4s)} = &-2 = -R_{12(4s)(4s-1)},\\
R_{12(4s-3)(4s-2)} =& -2 = -R_{12(4s-2)(4s-3)},\\
R_{13(4s)(4s-2)}=& -2 = -R_{13(4s-2)(4s)},\\
R_{13(4s-1)(4s-3)}=& 2=-R_{13(4s-3)(4s-1)},\\
R_{14(4s)(4s-3)} = &2 = -R_{14(4s-3)(4s)},\\
R_{14(4s-1)(4s-2)} = &2= -R_{14(4s-2)(4s-1)},
\end{align*}

and
$$
R_{1pAB} = 0, \qquad \text{otherwise}.
$$

Also,
\begin{align*}
d \eta_{1\alpha} &-  \eta_{1q} \wedge \eta_{q \alpha} - \eta_{1 \beta} \wedge \eta_{\beta \alpha} \\
&=d\eta_{\alpha} - 2\eta_q \wedge \eta_{q\alpha} - \eta_\beta \wedge \eta_{\beta \alpha}\\
& = -\eta_\alpha \wedge \eta_1 + e^t \, \omega_{q \alpha} \wedge
\eta_q,
\end{align*}
hence
\begin{align*}
R_{1 \alpha 1\alpha}&= -1,\\
R_{1(4s)(4s-1) 2} = &-1= -R_{1(4s-1)(4s)2},\\
R_{1(4s)(4s-2) 3} = &1 = -R_{1(4s-2)(4s)3},\\
R_{1(4s)(4s-3) 4} =& -1=-R_{1(4s-3)(4s)4},\\
R_{1(4s-1)(4s-3)3}=&-1=-R_{1(4s-3)(4s-1)3},\\
R_{1(4s-1)(4s-2)4} = &-1=-R_{1(4s-2)(4s-1)4},\\
R_{1(4s-2)(4s-3)2} = &-1=-R_{1(4s-3)(4s-2)2},
\end{align*}
and
$$
R_{1 \alpha AB} = 0 \qquad \text{otherwise}.
$$

Similarly,
\begin{align*}
d\eta_{pq} &- \eta_{p1} \wedge \eta_{1q} - \eta_{pr}\wedge\eta_{rq} - \eta_{p\beta} \wedge \eta_{\beta q}\\
&= d\omega_{pq} +4\eta_p \wedge \eta_q - \omega_{pr} \wedge \omega_{rq} - e^{2t}\omega_{p\beta} \wedge \omega_{\beta q}\\
&= \bar \Omega_{pq} + (1-e^{2t})\, \omega_{p \beta} \wedge
\omega_{\beta q} + 4 \eta_p \wedge \eta_q,
\end{align*}
where
$$
\bar \Omega_{pq} = \frac{1}{2}\bar R_{pqij} \omega_j \wedge \omega_i
$$
is the curvature form of $N_0.$ In particular, this implies that

\begin{equation}\label{eqn39}
R_{pqro} =\left\{ \aligned -4 + &e^{-4t} \,\bar R_{pqpq} \qquad \text{if} \qquad r=p \text{ and } o=q\\
4 + & e^{-4t}\,\bar R_{pqqp} \qquad \text{if} \qquad r=q \text{ and } o=p\\
&e^{-4t}\, \bar R_{pqro} \qquad \text{otherwise}, \endaligned
\right.
\end{equation}

\begin{align*}
R_{23(4s)(4s-3)} &= e^{-2t}\,\bar R_{23(4s)(4s-3)} - 2(e^{-2t} - 1),\\
R_{23(4s-1)(4s-2)} &= e^{-2t}\, \bar R_{23(4s-1)(4s-2)} - 2(e^{-2t} - 1),\\
R_{24(4s)(4s-2)} &= e^{-2t}\,\bar R_{24(4s)(4s-2)}  -2(e^{-2t} -1),\\
R_{24(4s-1)(4s-3)} &= e^{-2t} \, \bar R_{24(4s-1)(4s-3)} +2(e^{-2t} -1),\\
R_{34(4s)(4s-1)} &= e^{-2t}\, \bar R_{34(4s)(4s-1)} - 2(e^{-2t} -1),\\
R_{34(4s-2)(4s-3)} &= e^{-2t}\, \bar R_{34(4s-2)(4s-3)} - 2(e^{-2t}
-1),
\end{align*}
and
\begin{equation*}
    R_{pq\alpha \beta} = e^{-2t}\,\bar R_{pq \alpha \beta}, \qquad
\text{otherwise}.
\end{equation*}

We now continue with our curvature computation and consider
\begin{align*}
d\eta_{p\alpha} &- \eta_{p1} \wedge \eta_{1 \alpha} - \eta_{pq} \wedge \eta_{q\alpha} -\eta_{p\beta} \wedge \eta_{\beta \alpha}\\
&=d(e^t\, \omega_{p\alpha}) + 2 \eta_p \wedge \eta_\alpha - \omega_{pq} \wedge e^t\, \omega_{q \alpha} - e^t\, \omega_{p\beta} \wedge \eta_{\beta \alpha}\\
&= e^t\, \eta_1 \wedge \omega_{p \alpha} +\frac 12 e^t \bar
R_{p\alpha ij} \omega_j \wedge \omega_i+ 2 \eta_p \wedge \eta_\alpha
+ e^t\,\omega_{p \beta} \wedge (\omega_{\beta \alpha} -  \eta_{\beta
\alpha}),
\end{align*}
where $\bar R_{p \alpha ij}$ is the curvature tensor of $N_0.$ Using
$(\ref{eqn28}-\ref{eqn30})$ and $(\ref{eqn35}-\ref{eqn38})$, we have
\begin{align*}
 \frac12 &R_{2(4s) AB} \eta_B \wedge \eta_A\\
&= \eta_1 \wedge \eta_{(4s-1)} +\frac 12 e^t \bar R_{2(4s) ij}\,\omega_j\wedge \omega_i -2 \eta_{(4s)} \wedge \eta_2 + (1-e^{-2t})\, \eta_{(4s)} \wedge \eta_2 \\
&\qquad + (1-e^{-2t})\,\eta_{(4s-3)}\wedge \eta_3 + (1-e^{-2t})\, \eta_{(4s-2)} \wedge \eta_4\\
&= \eta_1 \wedge \eta_{(4s-1)}+\frac 12 e^t \bar R_{2(4s) ij}\,\omega_j\wedge \omega_i - (1+e^{-2t})\, \eta_{(4s)} \wedge \eta_2\\
&\qquad + (1-e^{-2t})\,\eta_{(4s-3)}\wedge \eta_3 + (1-e^{-2t})\,
\eta_{(4s-2)} \wedge \eta_4.
 \end{align*}

\begin{align*}
\frac12 &R_{2(4s-1) AB} \eta_B \wedge \eta_A\\
&=  -\eta_1 \wedge \eta_{(4s)} + \frac 12 e^t \bar R_{2(4s-1) ij}\,\omega_j\wedge \omega_i  - (1+e^{-2t})\,\eta_{(4s-1)}\wedge \eta_2 \\
&\qquad + (1-e^{-2t}) \eta_{(4s-2)} \wedge \eta_3 - (1-e^{-2t}) \,
\eta_{(4s-3)} \wedge \eta_4.
\end{align*}

\begin{align*}
\frac12 &R_{2(4s-2) AB} \eta_B \wedge \eta_A\\
&= \eta_1 \wedge \eta_{4s-3} +\frac 12 e^t \bar R_{2(4s-2) ij}\,\omega_j\wedge \omega_i  -(1+e^{-2t})\, \eta_{(4s-2)} \wedge \eta_2 \\
&\qquad + (1-e^{-2t})\,\eta_{(4s-1)}\wedge \eta_3 - (1-e^{-2t})\,
\eta_{(4s)} \wedge \eta_4.
\end{align*}

\begin{align*}
\frac12 & R_{2(4s-3) AB} \eta_B \wedge \eta_A\\
&= -\eta_1 \wedge \eta_{4s-2} +\frac 12 e^t \bar R_{2(4s-3) ij}\,\omega_j\wedge \omega_i - (1+e^{-2t})\, \eta_{(4s-3)} \wedge \eta_2\\
&\qquad - (1-e^{-2t})\,\eta_{(4s)}\wedge \eta_3 + (1-e^{-2t}) \,
\eta_{(4s-1)} \wedge \eta_4.
 \end{align*}

 Similar formulas for the curvature tensors of the form
 $ R_{3\alpha AB}$ and $R_{4\alpha AB}$.

Continuing with our computation of the second structural equations
using $(\ref{eqn28} -\ref{eqn30})$, we have
\begin{align}\label{eqn40}
d\eta_{(4s-1) (4s)} &- \eta_{(4s-1) 1}\wedge \eta_{1(4s)} - \eta_{(4s-1) q} \wedge \eta_{q (4s)} - \eta_{(4s-1) \beta} \wedge \eta_{\beta (4s)}\nonumber\\
&=d\omega_{(4s-1)(4s)} + 2e^{-2t}\, \eta_1 \wedge \eta_2  + (1-e^{-2t}) \,  d\eta_2 \nonumber\\
&\qquad+ \eta_{(4s-1)} \wedge \eta_{(4s)}- e^{2t}\, \omega_{(4s-1)q}\wedge \omega_{q(4s)}  \nonumber\\
&\qquad -(\omega_{(4s-1)(4s-2)} -(1-e^{-2t}) \, \eta_4) \wedge (\omega_{(4s-2) (4s)} - (1-e^{-2t})\, \eta_3)\nonumber\\
&\qquad  -(\omega_{(4s-1)(4s-3)}  - (1-e^{-2t})\, \eta_3) \wedge (\omega_{(4s-3)(4s)} + (1-e^{-2t})\, \eta_4)\nonumber\\
&= \frac12  \bar R_{(4s-1)(4s)ij}\, \omega_j \wedge \omega_i +(1-e^{2t})\, \omega_{(4s-1)q} \wedge \omega_{q(4s)} + 2 \eta_1 \wedge \eta_2  \nonumber\\
&\qquad+(1-e^{-2t})\,\omega_{2q} \wedge \eta_q + e^t(1-e^{-2t})\, \omega_{2\beta} \wedge \eta_\beta -\eta_{(4s)} \wedge \eta_{(4s-1)} \nonumber\\
&\qquad+ (1-e^{-2t})\, \omega_{(4s-1)(4s-2)} \wedge \eta_{3} + (1-e^{-2t})\, \eta_{4} \wedge \omega_{(4s-2)(4s)}  \nonumber\\
&\qquad+(1-e^{-2t}) \,\eta_3 \wedge \omega_{(4s-3)(4s)}  -(1 - e^{-2t})\, \omega_{(4s-1)(4s-3)}\wedge \eta_4  \nonumber\\
&\qquad +2(1-e^{-2t})^2 \, \eta_3\wedge \eta_4\nonumber\\
&= \frac 12 \bar R_{(4s-1)(4s)ij}\, \omega_j \wedge \omega_i +(2-e^{-2t})\, \eta_{(4s-1)} \wedge \eta_{(4s)} -2(1-e^{-2t})\, \eta_{(4s-3)} \wedge \eta_{(4s-2)}\nonumber\\
&\qquad + 2\eta_1 \wedge \eta_2+(1-e^{-2t})\, \omega_{2q} \wedge \eta_q  +2(1-e^{-2t})\,\eta_{(4r-3)} \wedge \eta_{(4r-2)} \nonumber\\
&\qquad + 2(1-e^{-2t})\,\eta_{(4r-1)} \wedge \eta_{(4r)}   + (1-e^{-2t})\,(\omega_{(4s-1)(4s-2)} - \omega_{(4s-3)(4s)}) \wedge \eta_3  \nonumber\\
&\qquad- (1-e^{-2t})\,(\omega_{(4s-2)(4s)} + \omega_{(4s-1)(4s-3)})
\wedge \eta_4+ 2(1-e^{-2t})^2\, \eta_3 \wedge \eta_4.
\end{align}
Note that $(\ref{eqn37})$ asserts that
\begin{align*} (1-e^{-2t})\, \omega_{2q} \wedge \eta_q
&= (1-e^{-2t})\,(-\omega_{14} \wedge \eta_3 + c \wedge \eta_3 + \omega_{13} \wedge \eta_4 - b \wedge \eta_4)\\
&= (1-e^{-2t})\,(4e^{-2t}\,\eta_3 \wedge \eta_4 +c \wedge \eta_3 - b
\eta_4),
\end{align*}
\begin{align*}
(1-e^{-2t})&\,(\omega_{(4s-1)(4s-2)} - \omega_{(4s-3)(4s)}) \wedge \eta_3\\
&=- (1-e^{-2t})\,c\wedge \eta_3,
\end{align*}
and
\begin{align*}
- (1-e^{-2t})&\,(\omega_{(4s-2)(4s)} + \omega_{(4s-1)(4s-3)}) \wedge \eta_4\\
&= (1-e^{-2t})\,b \wedge \eta_4.
\end{align*}
Hence substituting into $(\ref{eqn40})$, we obtain
\begin{align*}
\frac 12 &R_{(4s-1)(4s)AB}\,\eta_B \wedge \eta_A\\
&=\frac 12 \bar R_{(4s-1)(4s)ij}\, \omega_j \wedge \omega_i +(2-e^{-2t})\, \eta_{(4s-1)} \wedge \eta_{(4s)} -2(1-e^{-2t})\, \eta_{(4s-3)} \wedge \eta_{(4s-2)}+ 2\eta_1 \wedge \eta_2\\
&\qquad +2(1-e^{-2t})\,\eta_{(4r-3)} \wedge \eta_{(4r-2)} + 2(1-e^{-2t})\,\eta_{(4r-1)} \wedge \eta_{(4r)} +2(1-e^{-4t})\,\eta_3 \wedge \eta_4\\
&= 2\eta_1 \wedge \eta_2 + \frac 12 \bar R_{(4s-1)(4s)pq}\,e^{-4t}\,\eta_q \wedge \eta_p +2(1-e^{-4t})\,\eta_3 \wedge \eta_4 + \bar R_{(4s-1)(4s)p\alpha}\,e^{-3t}\,\eta_\alpha \wedge \eta_p \\
&\qquad + \frac 12 \bar R_{(4s-1)(4s)\alpha \beta}\,e^{-2t}\,\eta_\beta \wedge \eta_\alpha  +(2-e^{-2t})\, \eta_{(4s-1)} \wedge \eta_{(4s)} -2(1-e^{-2t})\, \eta_{(4s-3)} \wedge \eta_{(4s-2)}\\
& \qquad +2(1-e^{-2t})\,\eta_{(4r-3)} \wedge \eta_{(4r-2)} +
2(1-e^{-2t})\,\eta_{(4r-1)} \wedge \eta_{(4r)}.
\end{align*}

A similar computation yields the curvature tensor of the form
$R_{(4s-1)(4s-2)AB},$ $R_{(4s-1)(4s-3)AB},$  $R_{(4s-2)(4s-3)AB},$
$R_{(4s-2)(4s)AB},$ and $R_{(4s-3)(4s)AB}.$ It remains to compute
\begin{align}\label{eqn41}
\frac 12&R_{(4s-3)(4r)AB}\, \eta_B \wedge \eta_A\nonumber\\
&= d\eta_{(4s-3)(4r)} - \eta_{(4s-3)1}\wedge \eta_{1(4r)} - \eta_{(4s-3)q} \wedge \eta_{q(4r)} -\eta_{(4s-3)\beta} \wedge \eta_{\beta(4r)}\nonumber\\
&=d \omega_{(4s-3)(4r)} + \eta_{(4s-3)} \wedge \eta_{(4r)} - e^{2t}\, \omega_{(4s-3)q} \wedge \omega_{q (4r)} - \eta_{(4s-3) \beta} \wedge \eta_{\beta (4r)}\nonumber\\
&= \frac 12 \bar R_{(4s-3)(4r)ij}\, \omega_j \wedge \omega_i +(1-e^{2t})\, \omega_{(4s-3)q} \wedge \omega_{q(4r)} - (1-e^{-2t})\,(\eta_4 \wedge \omega_{(4s)(4r)} \nonumber\\
&\qquad+ \eta_3 \wedge \omega_{(4s-1)(4r)} ) -(1-e^{-2t})\,( \eta_2 \wedge \omega_{4s-2)(4r)} + \omega_{(4s-3)(4r-1)} \wedge \eta_2 \nonumber\\
& \qquad- \omega_{(4s-3)(4r-2)} \wedge \eta_3 +
\omega_{(4s-3)(4r-3)} \wedge \eta_4)  + \eta_{(4s-3)} \wedge
\eta_{(4r)}.
\end{align}
Using $(\ref{eqn27} -\ref{eqn29})$, we can write
$$
\omega_{(4s-3)q} \wedge \omega_{q(4r)} =-\eta_(4s-2) \wedge
\eta_{(4r-1)} + \eta_{(4s-1)} \wedge \eta_{(4r-2)} - \eta_{(4s)}
\wedge \eta_{(4r-3)}).
$$
Also using $(\ref{eqn27})$ asserts that

\begin{align*}
\omega_{(4s-3)(4r-1)} &= \omega_{(4s-2)(4r)},\\
\omega_{(4s-3)(4r-2)} &= -\omega_{((4s-1)(4r)}\\
\omega_{((4s-3)(4r-3)} &= \omega_{((4s)(4r)}.
\end{align*}

Hence $(\ref{eqn41})$ becomes

\begin{align*}
\frac 12&R_{(4s-3)(4r)AB}\, \eta_B \wedge \eta_A\\
&=\frac 12 \bar R_{(4s-3)(4r)pq}\,e^{-4t}\, \eta_q \wedge \eta_p + \bar R_{(4s-3)(4r)p\beta}\, e^{-3t}\,\eta_\beta \wedge \eta_p + \frac 12 \bar R_{(4s-3)(4r)\alpha \beta}\,e^{-2t}\, \eta_\beta \wedge \eta_\alpha \\
& \qquad -(1-e^{2t})\,\eta_{(4s-2)} \wedge \eta_{(4r-1)}
+(1-e^{-2t})\, \eta_{(4s-1)} \wedge \eta_{(4r-2)} \nonumber\\
& \qquad - (1-e^{-2t})\,\eta_{(4s)} \wedge \eta_{(4r-3)} +
\eta_{(4s-3)} \wedge \eta_{4r}.
\end{align*}

So we have determined all curvature tensors of $M$.  Note that the
quaternionic curvatures satisfy
\begin{align*}
\K(e_1, e_2) + \K(e_1, e_3) + \K(e_1, e_4) &= -12\\
\K(e_2, e_1) + \K(e_2, e_3) + \K(e_2, e_4) &= -12 + e^{-2t}\,(\K^N(e_2, e_3) + \K^N(e_2, e_4))\\
\K(e_3, e_1) + \K(e_3, e_2) + \K(e_3, e_4) &= -12 + e^{-2t}\,(\K^N(e_3, e_2) + \K^N(e_3, e_4))\\
\K(e_4, e_1) + \K(e_4, e_2) + \K(e_4, e_3) &= -12 + e^{-2t}\,(\bar
\K(e_4, e_2) + \K^N(e_4, e_3)).
\end{align*}
In particular, this implies that
$$
\K^N(e_2, e_3) = \K^N(e_2, e_4) = \K^N(e_3, e_4) = 0.
$$
Also, for $2 \le p \le 4,$ we have
\begin{align*}
\sum_{i=0}^3 \K(e_1, e_{(4s-i)}) &= -4\\
\sum_{i=0}^3 \K(e_p, e_{(4s-i)} & = -4 + e^{-2t}\,(\sum_{i=0}^3 \bar
\K(e_p, e_{(4s-i)}) - 4)
\end{align*}
implying
$$
\sum_{i=0}^3 \K^N(e_p, e_{(4s-i)}) = 4.
$$
We also have
$$
\sum_{i=1}^3 \K(e_{(4s)}, e_{(4s-1)}) = -12 + e^{-2t}\,(\sum_{i=1}^3
\K^N(e_{(4s)}, e_{(4s-i)}) +9)
$$
implying
$$
\sum_{i=1}^3 \K^N(e_{(4s)}, e_{(4s-i)}) = -9.
$$
Lastly,
$$
\sum_{i=0}^3 \K(e_{(4s)}, e_{(4r-i)}) = -4 + e^{-2t}\,\sum_{i=0}^3
\K^N(e_{(4s)}, e_{(4r-i)})
$$
implying
$$
\sum_{i=0}^3 \K^N(e_{(4s)}, e_{(4r-i)}) = 0.
$$

The above computation determined the whole curvature tensor for $M$
and $N_0.$  In particular, if $M$ has bounded curvature, then from
the  formulas about the components of curvature tensors of $M$, all
curvature components are determined as those of $\mathbb{QH}^n$. So
it must be covered by $ {\mathbb{QH}}^n.$ \QED

\noindent Shengli Kong\newline Department of Mathematics\newline
University of California, Irvine\newline Irvine, CA92697-3875, USA
\newline email: skong@math.uci.edu
\medskip

\noindent Peter Li\newline Department of Mathematics\newline University of California, Irvine\newline Irvine, CA92697-3875, USA
\newline email:pli@math.uci.edu

\medskip
\noindent Detang Zhou\newline Departamento de Geometria\newline
Insitituto de Matematica\newline Universidade Federal Fluminense-
UFF\newline Centro, Niter\'{o}i, RJ 24020-140, Brazil \newline
email: zhou@impa.br

\enddocument